ANNALES
DE L'INSTITUT
HENRI
POINCARÉ
PROBABILITÉS
ET STATISTIQUES

# A local limit theorem with speed of convergence for Euclidean algorithms and diophantine costs

## Viviane Baladi  and Aïcha Hachemi

*CNRS-UMR 7586, Institut de Mathématiques Jussieu, Paris, France.*
*E-mail: baladi@math.jussieu.fr; hachemi@math.jussieu.fr*

**Abstract.** For large $N$, we consider the ordinary continued fraction of $x = p/q$ with $1 \le p \le q \le N$, or, equivalently, Euclid's gcd algorithm for two integers $1 \le p \le q \le N$, putting the uniform distribution on the set of $p$ and $q$s. We study the distribution of the total cost of execution of the algorithm for an additive cost function $c$ on the set $\mathbb{Z}_+^*$ of possible digits, asymptotically for $N \to \infty$. If $c$ is nonlattice and satisfies mild growth conditions, the local limit theorem was proved previously by the second named author. Introducing diophantine conditions on the cost, we are able to control the speed of convergence in the local limit theorem. We use previous estimates of the first author and Vallée, and we adapt to our setting bounds of Dolgopyat and Melbourne on transfer operators. Our diophantine condition is generic (with respect to Lebesgue measure). For smooth enough observables (depending on the diophantine condition) we attain the optimal speed.

**Résumé.** Nous considérons la fraction continue ordinaire de $x = p/q$ pour $1 \le p \le q \le N$, ou, de manière équivalente, l'algorithme de pgcd d'Euclide pour deux entiers $1 \le p \le q \le N$, avec $N$ grand et $p$ et $q$ distribués uniformément. Nous étudions la distribution du coût total de l'exécution de l'algorithme pour un coût additif $c$ sur l'ensemble $\mathbb{Z}_+^*$ des "digits" possibles, lorsque $N$ tend vers l'infini. Le théorème de la limite locale a été démontré par le deuxième auteur si $c$ est non réseau et satisfait une condition de croissance modérée. En imposant une condition diophantienne sur le coût, nous parvenons à contrôler la vitesse de convergence dans ce théorème de la limite locale. Pour cela nous utilisons des estimées obtenues par le premier auteur et Vallée, et nous adaptons à notre problème des bornes de Dolgopyat et Melbourne sur les opérateurs de transfert. Notre condition diophantienne est générique (par rapport à la mesure de Lebesgue). Pour des observables assez régulières (par rapport à la condition diophantienne), nous obtenons la vitesse optimale.

# 1. Introduction and statement of results

Every rational $x \in {]0,1]}$ admits a finite continued fraction expansion

$$x = \cfrac{1}{m_1 + 1/(m_2 + 1/(\cdots + 1/m_P))}, \quad m_j = m_j(x) \in \mathbb{Z}_+^*, P = P(x) \in \mathbb{Z}_+^*. \tag{1}$$

Continued fraction expansions can be viewed as trajectories of the Gauss map

$$T : (0,1] \to [0,1], \qquad T(x) := \frac{1}{x} - \left[\frac{1}{x}\right].$$

(Here, $[y]$ is the integer part of $y \in \mathbb{R}_+^*$.) Indeed, if $x \neq 0$ is rational, then $T^P(x) = 0$ for some $P = P(x) \geq 1$, which is the depth of the continued fraction, and the digits $m_j = m_j(x) \in \mathbb{Z}_+^*$ appearing in (1) are just

$$m_j(x) = \left[\frac{1}{T^{j-1}(x)}\right], \quad 1 \leq j \leq P(x).$$

Clearly, this is equivalent to execution of Euclid's gcd algorithm: for two integers $1 \leq p \leq q$, write $q_1 = q$, $p_1 = p$ and $q_1 = m_1 p_1 + r_1$ with $m_1 = m_1(p/q) \in \mathbb{Z}_+^*$ and a remainder $r_1 \in \mathbb{Z}_+$ so that $r_1 < p_1$. If $r_1 = 0$ we are done, and $p = \gcd(p,q)$, with $P(p/q) = 1$. If $r_1 \neq 0$, set $p_2 = r_1$ and $q_2 = p_1$, and iterate this procedure until the remainder $r_P$ vanishes for some $P = P(p/q) \geq 2$. Then $p_P = \gcd(p,q)$, and $m_j = m_j(p/q)$ for $1 \leq j \leq P$. Note that $m_{P(p/q)} = 1$ if and only if $p = q$.

We shall call *cost* any (nonidentically zero) function $c : \mathbb{Z}_+^* \to \mathbb{R}$. Given such $c$, we associate to each rational $x = p/q \in (0,1]$ the following *total cost*:

$$C(x) = \sum_{j=1}^{P(x)} c(m_j(x)). \tag{2}$$

(Note that if $c \equiv 1$ then the total cost is just the depth $P(x)$ of the continued fraction.) Our goal is to describe the probabilistic behaviour of the total cost associated to the ordinary (Gauss) continued fraction (1) (and some of its "fast" variants). Before stating our results, we explain our probabilistic setting, and recall previous works.

Consider $\tilde{\Omega} := \{(p,q) \in (\mathbb{Z}_+^*)^2, \frac{p}{q} \in (0,1]\}$, and $\Omega := \{(p,q) \in \tilde{\Omega} \mid \gcd(p,q) = 1\}$, and endow the sets $\tilde{\Omega}_N := \{(p,q) \in \tilde{\Omega} \mid q \leq N\}$, and $\Omega_N := \{(p,q) \in \Omega \mid q \leq N\}$ with uniform probabilities $\tilde{\mathbb{P}}_N$ and $\mathbb{P}_N$, respectively. We shall state our results for $\mathbb{P}_N$, but, as observed e.g. in [1] (see (2.18)), they also hold for $\tilde{\mathbb{P}}_N$. Note that if $(p,q) \in \Omega_N$ we can write $C(p,q)$ and $P(p,q)$ instead of $C(p/q)$ and $P(p/q)$. As usual, the expectation $\mathbb{E}_N(C)$ denotes $\sum_{(p,q) \in \Omega_N} \mathbb{P}_N((p,q))C(p,q)$ and the variance $\mathbb{V}_N(C)$ is $\mathbb{E}_N(C^2) - (\mathbb{E}_N(C))^2$.

We shall use the following conditions on a cost function $c$:

$$\begin{cases} c \text{ is of moderate growth if there exists } \nu > 0 \text{ so that } \sum_{m \in \mathbb{Z}_+^*} e^{\nu c(m)} m^{-2+\nu} < \infty, \\ c \text{ has strong moments up to order } k \geq 1 \text{ if there exists } \nu > 0 \text{ so that } \sum_{m \in \mathbb{Z}_+^*} c(m)^k m^{-2+\nu} < \infty. \end{cases}$$

Of course, if $c$ is of moderate growth then it has strong moments up to arbitrary order $k$.

**Remark 1.1.** *If $c(m) = \mathrm{O}(\log(m))$ then $c$ is of moderate growth, while if $c(m) = \mathrm{O}(m^{-\nu'+1/k})$ for some $k \geq 1$ and some $\nu' > 0$ then $c$ has strong moments up to order $k$. The terminology comes from the fact that $T$ has a unique absolutely continuous invariant measure $\mu_1$ on $[0,1]$, with a positive analytic density $f_1$ (see also Section 2), so that, writing $c_*(x) = c(m)$ if $x \in (1/(m+1), 1/m]$, we have that $\int (c_*(x))^\ell \,\mathrm{d}\mu_1(x)$ is well-defined for positive integers $\ell \leq k$, i.e. $c_*$ has moments up to order $k$ if $c$ has strong moments up to order $k$. (The converse is not true, whence the terminology "strong" moments. See [1], and Lemma 3.3 for the need to use $\nu > 0$.)*

Introducing a dynamical approach (centered around the one-dimensional map $T$) and using Ruelle-type transfer operators $\mathbf{H}_s = \mathbf{H}_{s,0}$ to study this problem (see Section 2 for a definition of the transfer operators), Brigitte Vallée obtained in a series of papers (see e.g. [21] for references) precise results for the asymptotics of the expectation $\mathbb{E}_N(C)$ and other moments. For example, if $c$ is of moderate growth, there is $\mu(c) \in \mathbb{R}$ (with $\mu(c) \neq 0$ if $\int c_*(x)\,d\mu_1(x) \neq 0$) so that

$$\lim_{N \to \infty} \frac{\mathbb{E}_N(C)}{\log N} = \mu(c). \tag{3}$$

We refer to the recent work [1] for more information, a historical discussion including references to the work of Heilbronn and Dixon and previous work of Vallée. Among other things it was proved in [1] that if $c$ is of moderate growth then there exists $\delta(c) \in \mathbb{R}_+^*$ so that

$$\lim_{N \to \infty} \frac{\mathbb{V}_N(C)}{\log N} = (\delta(c))^2. \tag{4}$$

**Remark 1.2.** *To get (3) and (4) it suffices to assume that $c$ has strong moments up to order 2. See Lemma 3.2.*

The article [1] also contains the following *central limit theorem* [1], Theorem 3: for each cost $c$ of moderate growth, there is $M_1(c) \geq 1$ so that for any integer $N \geq 1$, and any $y \in \mathbb{R}$:

$$\left| \mathbb{P}_N\left( \frac{C(p,q) - \mu(c)\log N}{\delta(c)\sqrt{\log N}} \leq y \right) - \frac{1}{\sqrt{2\pi}} \int_{-\infty}^{y} e^{-x^2/2}\,dx \right| \leq \frac{M_1(c)}{\sqrt{\log N}}, \tag{5}$$

where $\mu(c) \in \mathbb{R}$ and $\delta(c) > 0$ are the same as in (3). The speed of convergence $(\log N)^{-1/2}$ in the above central limit theorem is optimal, as is clear from the saddle point-argument in the proof. This speed is the equivalent in our setting of the speed of convergence in the central limit theorem for independent identically distributed random variables [10]. (Hensley [14] obtained a central limit theorem for $c \equiv 1$, more than a decade before [1], but with a $O((\log N)^{-1/24})$ bound on the rate of convergence.)

The basic tool to obtain all limit theorems mentioned here is the transfer operator. The transfer operator allows to implement the characteristic function method, via Levy generating functions. However, there is a key difference between *discrete* and *continuous* problems (see also [1], Section 2.3, for an illustration). Discrete and continuous does not refer to time here, but to the probabilities: discrete means that we are performing weighted sums over finite sets of increasing size (just like in this paper). In the continuous setting, such as in the pioneering work of Guivarc'h–Le Jan [12], in a geometric context which also involves the Gauss map, a spectral gap argument à la Nagaev is invoked. However, the transfer operators which appear for discrete problems involve not only the parameter $i\tau$ or $w$ from the characteristic function, but also another parameter, $s$, which ranges in a half-plane containing the pole. This other parameter comes from a Dirichlet series in the present paper and in [1], it can also be viewed as the parameter of an $L$-function in other settings, see e.g. [19]. Dealing with $s$ of large imaginary parts requires the use of fundamental bounds first proved by Dolgopyat [7, 8] in the context of hyperbolic flows.

Specifically, in order to prove (5), methods adapted[1] from Dolgopyat [8], were used to get bounds on the resolvent $(\mathrm{Id} - \mathbf{H}_{s,w})^{-1}$ of the transfer operator $\mathbf{H}_{s,w}$, with $(s,w) \neq (1,0)$, when the complex parameter $s$ varies in a half-plane containing 1, but the complex parameter $w$ is close to zero. We refer to [1] (in particular Theorem 2 there) for more details.

We next discuss local limit theorems. A cost function $c$ is called *lattice*, if there exists $(L, L_0) \in \mathbb{R}_+^* \times \mathbb{R}_+$, with $L_0/L \in [0,1)$ so that $(c - L_0)/L$ is integer-valued. If $c$ is lattice but not constant, the largest possible $L$ is called the *span* of $c$ and the corresponding $L_0$ is called the *shift* of $c$. If $c$ is constant we take span $L = |c|$ and shift $L_0 = 0$.[2] A cost function is called *nonlattice* if it is not lattice.

---

[1] The key difficulty is that the symbolic alphabet in [8] is finite, while the Gauss map has infinitely many branches.

[2] The definition of lattice stated in [1] and [13] should be replaced there by this one.

If the cost is lattice with $L_0 = 0$ and enjoys moderate growth, then the following *local limit theorem* ([1], Theorem 4) holds: for $x \in \mathbb{R}$ and $N \in \mathbb{Z}_+^*$, put

$$\mathbf{Q}(x, N) = \mu(c) \log N + \delta(c) x \sqrt{\log N}, \tag{6}$$

then,[3] there is $M_2(c) \geq 1$ so that for every $x \in \mathbb{R}$ and all integers $N \geq 1$

$$\left| \sqrt{\log N} \cdot \mathbb{P}_N \left( (C(p,q) - \mathbf{Q}(x, N)) \in \left( -\frac{L}{2}, \frac{L}{2} \right] \right) - \frac{L e^{-x^2/2}}{\delta(c)\sqrt{2\pi}} \right| \leq \frac{M_2(c)}{\sqrt{\log N}}. \tag{7}$$

(See [1], Section 5.4, for the case $L_0 \neq 0$.) Again, the constants $\mu(c) \in \mathbb{R}$ and $\delta(c) > 0$ are the same as in (3), and the speed of convergence is optimal. The proof uses operators $\mathbf{H}_{s, i\tau}$ where the complex parameter $s$ varies in a half-plane, and the real parameter $\tau$ lies in the bounded interval $[-\pi/L, \pi/L]$.

Very recently, using Breiman's method (also known as Stone's trick, see [11] for a recent application of this method to nonuniformly hyperbolic dynamics) to handle noncompactness issues, the second author of the present paper obtained [13], Théorème 3, a local limit theorem: for every nonlattice[4] cost function $c$ of moderate growth, for each compact interval $J$, we have, writing $|J|$ for the length of $J$, and for $\mu(c) \in \mathbb{R}$ and $\delta(c) > 0$ as in (3):

$$\lim_{N \to \infty} \sup_{x \in \mathbb{R}} \left| \sqrt{\log N} \cdot \mathbb{P}_N ((C(p,q) - \mathbf{Q}(x, N)) \in J) - |J| \frac{e^{-x^2/2}}{\delta(c)\sqrt{2\pi}} \right| = 0. \tag{8}$$

We would like to point out that the compact interval $J$ in (8) is arbitrary, while in (7) only the interval $(-L/2, L/2]$ is allowed. Roughly speaking, the local limit theorem in the lattice case can be viewed as the analogue of a result for a Poincaré section of a flow. (And it is well-known that obtaining results for flows is usually much more difficult, sometimes requiring additional conditions, than proving the corresponding results for their Poincaré maps.)

**Remark 1.3.** *Lemma 3.2 and the arguments in Section 4 show that the moderate growth assumption for (5), (7) and (8) can be replaced by the requirement that c has strong moments up to order 3.*

The purpose of the present article is to obtain a local limit theorem with control of the speed of convergence in the nonlattice case.

Our proof involves transfer operators $\mathbf{H}_{s, i\tau}$, with $s$ in a half-plane, and $\tau \in \mathbb{R}$. When $\tau$ is confined to any compact set, or when $|\Im s|$ is large enough, we can exploit the bounds of [1] based on Dolgopyat's work [8] (see Lemmas 3.3 and 3.4). For large $\tau$, we must adapt other estimates of Dolgopyat [7], introduced to study decay of correlations for Axiom A flows (see also Naud [16] for a reader-friendly account). These estimates of Dolgopyat [7] require a diophantine condition.[5] In view of proving decay of correlations for nonuniformly hyperbolic flows, Melbourne [15] recently generalised the estimates of [7] to infinite alphabets. However, the parameters $(\tilde{s}, z)$ of Melbourne's operators $R_{\tilde{s}, z}$ [15], Section 3.3, are not of the same nature as our parameters $(s, i\tau)$: while the real and imaginary parts of $\tilde{s}$ correspond to $\Re \tilde{s} = \Re s$ and $\Im \tilde{s} = \tau$, respectively, the imaginary part of $z$ lies in $[0, 2\pi)$, because it arises from locally constant integer return times. The parameter $\Im z$ is thus a periodic parameter, and the return times are "mute" in a reformulation of Melbourne's diophantine condition ([15], Proof of Corollary 2.4). In our setting, unbounded values of $\Im s$ are handled by the arguments from [1], as mentioned above, and we are left to deal with $|\Im s|$ in a compact set. The parameter $\Im s$ is thus "artificially" bounded, but is not intrinsically periodic. Because of these differences (note also that Melbourne works with the dynamical distance in an abstract Gibbs–Markov setting, while we require the Euclidean metric), we carry out the modified bounds in detail in Section 2. Since the weights $|h'|$ are not locally constant integers, we cannot eliminate them from our diophantine condition (see also Remark 1.4 and the proof of Lemma 1.5).

---

[3]The factor $L$ multiplying $e^{-x^2/2}$ was inadvertently omitted from [1], pp. 351 and 381.

[4]The assumption that $c$ is nonlattice has been inadvertently omitted from [13], Théorème 3: this assumption is necessary to allow arbitrarily large $L$ in [13], Lemma 2, as is clear e.g. from the proof of [1], Proposition 1(iii).

[5]See [2, 3] for diophantine conditions in a related, but different, probabilistic context, and [4] for a previous "nonlattice of order $p$" condition.

Since the statement of our sufficient diophantine condition on the cost is unpalatable (although very similar to Melbourne's [15], Theorem 2.3), we postpone it to Section 2, and we formulate here a stronger (but simpler) condition instead. We need more notation. The countable set of digits $m \in \mathbb{Z}_+^*$ is in bijection with the set $\mathcal{H}$ of inverse branches of $T$, through $m \mapsto (y \mapsto \frac{1}{y+m})$. We may thus view the cost function $c$ as a function on $\mathcal{H}$. For integer $p \geq 2$ and any subset $\mathcal{H}_0$ of $\mathcal{H}$, write $\mathcal{H}_0^1 = \mathcal{H}_0$ and $\mathcal{H}_0^p = \{h \circ \tilde{h} \mid h \in \mathcal{H}_0, \tilde{h} \in \mathcal{H}_0^{p-1}\}$. Note that $x = h(x)$ for $x \in \mathcal{H}^p$ means that $x = T^p(x)$, i.e. $x$ is periodic. The minimal such $p$ is called the period of $x$. Then, we extend $c$ to a function on $\mathcal{H}^p$ by setting, for $h = h_p \circ \cdots \circ h_2 \circ h_1 \in \mathcal{H}^p$,

$$c(h) = \sum_{j=1}^{p} c(h_j). \tag{9}$$

Recall that a vector $x \in \mathbb{R}^d$ for $d \geq 1$ is *diophantine* of exponent $\eta_0 \geq d$ if there exists $M > 0$ so that for each $(q_1, \ldots, q_d) \in \mathbb{Z}^d \setminus \{0\}$

$$\inf_{p \in \mathbb{Z}} \left| p - \sum_{k=1}^{d} x_k q_k \right| \geq \frac{M}{(\max_k |q_k|)^{\eta_0}}.$$

For each $\eta_0 > d$, the set of diophantine vectors of exponent $\eta_0$ has full Lebesgue measure in $\mathbb{R}^d$ (see e.g. [5]).

**Definition.** *Let $\eta > 2$. The cost $c$ is strongly diophantine of exponent $\eta$ if there exist $\eta > \eta_0 \geq 2$ and four periodic points $x_j \in (0,1)$, $j = 1, 2, 3, 4$ for $T$, of respective minimal periods $p_j \geq 1$, with $h_{i_j}(x_j) = x_j$ for $h_{i_j} \in \mathcal{H}^{p_j}$, and with pairwise disjoint orbits, so that the following holds: setting $c_j = c(h_{i_j})$, $a_j = \log|h'_{i_j}|$,*

$$L_{1j} = p_j c_1 - p_1 c_j, \qquad \hat{L}_{1j} = p_j a_1 - p_1 a_j, \quad j = 2, 3, 4,$$

*then $L_{13} \neq 0$, $\hat{L}_{12} \neq 0$, with $L_{12}/L_{13}$ diophantine of exponent $\eta_0$, and, defining*

$$\tilde{L}_{jk} = L_{1j}\hat{L}_{1k} - \hat{L}_{1j}L_{1k}, \quad j \neq k \in \{2, 3, 4\},$$

*we have $\tilde{L}_{23} \neq 0$, and the pair*

$$\left( \frac{\tilde{L}_{43}}{\tilde{L}_{23}}, \frac{\tilde{L}_{42}}{\tilde{L}_{23}} \right)$$

*is diophantine of exponent $\eta_0$.*

**Remark 1.4.** *Our definition uses four periodic orbits, like in [15], with the "intertwining" of the $c_j$ and the $a_j$ due to the previously discussed fact that the $a_j$ are not integers. (See also Remark 1.6.) The condition $L_{13} \neq 0$ with $L_{12}/L_{13}$ diophantine of exponent $\eta_0$, involving three periodic orbits, is sufficient to ensure that $\|(\mathrm{Id} - \mathbf{H}_{\sigma,i\tau})^{-1}\|_{\mathrm{Lip}} \leq M_0 |\tau|^\alpha$ for some $\alpha > 0$, all $|\tau| \geq 2$ and appropriate $\sigma$ (see proof of Proposition 2.1), this is not enough since we need $\|(\mathrm{Id} - \mathbf{H}_{\sigma+it,i\tau})^{-1}\|_{\mathrm{Lip}} \leq M_0 |\tau|^\alpha$ for all $|t| \leq t_0$. In the simpler cases studied by Dolgopyat [7] and Naud [16] it is possible to formulate a (non intertwined) sufficient diophantine condition involving only two periodic orbits (because the alphabet is finite).*

Before stating our result, we discuss further the above definition. First, it is easy to see that a strongly diophantine cost is nonlattice. By Remark 1.8, it is generic. In Section 2, we shall give the definition of a diophantine cost of exponent $\eta$ and we shall prove:

**Lemma 1.5.** *If $c$ is strongly diophantine of exponent $\eta$, letting $\mathcal{H}_0 \subset \mathcal{H}$ be the smallest set so that $h_{i_j} \in \mathcal{H}_0^{p_j}$ for $j = 1, 2, 3, 4$, then $c$ is diophantine of exponent $\eta$ for $\mathcal{H}_0$.*

**Remark 1.6.** *Each $x_j$ in the strongly diophantine condition is just the quadratic number associated to the infinite repetition of a $p_j$-tuple $\vec{m} = (m_{1,j}, \ldots, m_{p_j,j})$ of positive integers. Also, $a_j$ coincides with $\log \prod_{\ell=0}^{p_j-1}(m_{\ell,j} + x_{j,\ell})^{-2}$, where the $x_{j,\ell} = T^\ell(x_j)$ are the quadratic numbers associated to the circular permutations of $\vec{m}$. Finally $c_j$ is the total cost associated to the rational number $u/v$ whose continued fraction has depth $p_j$ and is given by $\vec{m}$.*

**Remark 1.7.** *We give an example. If $c$ is such that there are four integers $m_j > 0$ so that, setting $a_j = \log(1 + (m_j^2 - m_j\sqrt{m_j^2 + 4})/2)$, we have*

$$c(m_2) \neq c(m_3), \qquad a_1 \neq a_2, \qquad (c(m_1) - c(m_2))(a_1 - a_3) \neq (c(m_1) - c(m_3))(a_1 - a_2),$$

*and both $\frac{c(m_1)-c(m_2)}{c(m_1)-c(m_3)}$ and the pair*

$$\left( \frac{(c(m_1) - c(m_4))(a_1 - a_3) - (c(m_1) - c(m_3))(a_1 - a_4)}{(c(m_1) - c(m_2))(a_1 - a_3) - (c(m_1) - c(m_3))(a_1 - a_2)}, \right.$$

$$\left. \frac{(c(m_1) - c(m_4))(a_1 - a_2) - (c(m_1) - c(m_2))(a_1 - a_4)}{(c(m_1) - c(m_2))(a_1 - a_3) - (c(m_1) - c(m_3))(a_1 - a_2)} \right)$$

*are diophantine of exponent $\eta_0 < \eta$, then $c$ is strongly diophantine of exponent $\eta$. (Just consider the case where all $p_j = 1$ in the definition.)*

**Remark 1.8.** *Fix $\eta > 2$. Choose four periodic points $x_j = h_{i_j}(x_j) \in (0, 1]$ of $T$, with pairwise disjoint orbits (and $h_{i_j} \in \mathcal{H}^{p_j}$). Then it is not difficult (see [15], Corollary 2.4) to show that for Lebesgue almost every $(c_1, c_2, c_3, c_4)$ in $\mathbb{R}_+^4$, any cost $c$ so that $c(h_{i_j}) = c_j$ is strongly diophantine of exponent $\eta$ (use Fubini). Therefore, the diophantine condition on the cost $c$ deserves to be called generic if $\eta > 2$.*

Recall (6). Our main result is the following local limit theorem with speed of convergence:

**Theorem 1.9.** *For any diophantine cost function $c$ of exponent $\eta$ and subset $\mathcal{H}_0$, with strong moments up to order 3:*

*There exists $\varepsilon \in (0, 1/2]$ so that for each compact interval $J \subset \mathbb{R}$ there exists a constant $M_J > 0$ so that for every $x \in \mathbb{R}$ and all integers $N \geq 1$:*

$$\left| \sqrt{\log N}\, \mathbb{P}_N((C(p,q) - \mathbf{Q}(x,N)) \in J) - |J| \frac{e^{-x^2/2}}{\delta(c)\sqrt{2\pi}} \right| \leq \frac{M_J}{(\log N)^\varepsilon}.$$

*There exists $r \geq 1$ so that for any compactly supported $\psi \in C^r(\mathbb{R})$, there exists a constant $M_\psi > 0$ so that for every $x \in \mathbb{R}$ and all integers $N \geq 1$:*

$$\left| \sqrt{\log N}\, \mathbb{E}_N(\psi(C(p,q) - \mathbf{Q}(x,N))) - \frac{e^{-x^2/2}}{\delta(c)\sqrt{2\pi}} \int \psi(y)\,\mathrm{d}y \right| \leq \frac{M_\psi}{\sqrt{\log N}}.$$

The second claim of the theorem says that we attain the optimal speed in the local limit theorem for smooth enough compactly supported observables $\psi$. If the cost is nonlattice but not diophantine, we expect that arbitrarily slow convergence can take place in the local limit theorem, in the spirit of [17, 18].[6]

---

[6]Note however that lower bounds are much less accessible in our setting and the methods in this paper do not provide such bounds for any examples of non-diophantine lattice costs.

**Remark 1.10.** *The statement of Proposition 2.1 and the proof of Theorem 1.9 easily imply that for any*

$$\alpha > \eta\left(2 + \frac{\log\sup_{\mathcal{H}_0}|h'|^{-1}}{\log(1/\rho)}\right)\left(1 + \frac{\log\sup_{\mathcal{H}_0}|h'|^{-1}}{\log(1/\rho)}\right),$$

*we may take*

$$\varepsilon < (2\alpha)^{-1}, \qquad r > \alpha + 1, \tag{10}$$

*but these conditions are probably not optimal. The proof of Theorem 1.9 gives a constant $K(c) = K(\eta, \mathcal{H}_0)$ so that, if $\psi$ is supported in an interval $J_\psi$,*

$$M_J \leq K(c)|J|^2, \qquad M_\psi \leq K(c)|J_\psi| \cdot \|\psi\|_{C^{[r+1]}}. \tag{11}$$

**Remark 1.11.** *If $c$ has strong moments up to order $k+1 \geq 4$, a little more work should yield finite Edgeworth expansions [10] (see also [2]) of order $k$ for compactly supported $\psi \in C^r(\mathbb{R})$. (The remainder term being $O((\log N)^{-k/2})$.) In this case, the condition on the differentiability $r$ of $\psi$ will depend not only on the diophantine exponent of $c$, but also on the desired order $k$ for the Edgeworth expansion.*

The paper is organised as follows: in Section 2, we adapt the estimates of Dolgopyat–Melbourne ([7, 15]) to our setting to get bounds (Proposition 2.1) for the norm of the resolvent $(\mathrm{Id} - \mathbf{H}_{\sigma+it,i\tau})^{-1}$ for large $|\tau|$, bounded $|t|$, and $\sigma > 1 - \delta(\tau)$ for small $\delta(\tau)$, under the diophantine assumption on $c$. Lemma 1.5 is also proved in Section 2. In Section 3, we first recall previous material from [1], in particular the connection between $(\mathrm{Id} - \mathbf{H}_{\sigma+it,i\tau})^{-1}$ and the moment generating function $\overline{\mathbb{E}}_N(\mathrm{e}^{i\tau C})$ of smoothened models (via the Perron formula and bivariate Dirichlet series), as well as estimates on $\overline{\mathbb{E}}_N(\mathrm{e}^{i\tau C})$ for bounded $|\tau|$. We then deduce from Proposition 2.1 our key estimate (Corollary 3.5), on $\overline{\mathbb{E}}_N(\mathrm{e}^{i\tau C})$ for large $|\tau|$. The proof of Theorem 1.9 is carried out in Section 4 by reducing to a study of $\int_{\mathbb{R}} \hat{\chi}_J(\tau)\mathrm{e}^{-i\tau \mathbf{Q}(x,N)}\overline{\mathbb{E}}_N(\mathrm{e}^{i\tau C})\,d\tau$, respectively $\int_{\mathbb{R}} \hat{\psi}(\tau)\mathrm{e}^{-i\tau \mathbf{Q}(x,N)}\overline{\mathbb{E}}_N(\mathrm{e}^{i\tau C})\,d\tau$ ($\hat{\phi}$ the Fourier transform of $\phi$), and decomposing the integral over $\tau \in \mathbb{R}$ into four domains, over which we apply the estimates from Section 3. In the Appendix, we describe two other (fast) continued fraction algorithms, the *centered algorithm* and the *odd algorithm,* for which Theorem 1.9 holds, with the same proof, since our arguments only use the fact that $T$ belongs to the good class from [1] and satisfies the condition UNI from [1].

We have already mentioned the difference between the bounds in Section 3 and those in [15]. With respect to [1] and [13], we can mention two innovations (besides our remark that moderate growth can be replaced by strong moments): first, Lemma 3.5 requires a specific smoothening, adapted to the weak bounds for large $|\tau|$ from Section 3; second, we have to regularise the characteristic function of the interval $J$ by convolution in Section 4.2 in order to control large $|\tau|$.

## 2. Dolgopyat–Melbourne estimates for $\mathbf{H}_{\sigma+it,i\tau}$

Let us first introduce some notation. Put $I = [0,1]$, and $\mathrm{Lip}(I) = \{u : I \to \mathbb{C} \mid \|u\|_{L^\infty} + \mathrm{Lip}(u) < \infty\}$, with $\mathrm{Lip}(u)$ the smallest Lipschitz constant of $u$. (If $u$ is not Lipschitz then we put $\mathrm{Lip}(u) = \infty$.) It is well known (see [1] for references) that there exists $\rho < 1$, $K \geq 1$ and $\widehat{K} \geq 1$ so that for all $m \in \mathbb{Z}_+^*$ and all $h \in \mathcal{H}^m$

$$\sup|h'| \leq K\rho^m, \qquad |h''(x)| \leq \widehat{K}|h'(x)|, \quad \forall x \in I. \tag{12}$$

In this section, we focus on $|t| \leq t_0$, for some fixed $t_0 > 0$, and $|\tau| \geq 2$, since other values of $t$ and $\tau$ are covered in previous works, as explained in the next section.

For $\tau \in \mathbb{R}$ and $s = \sigma + it$, with $t \in \mathbb{R}$, $\sigma \in \mathbb{R}$ with $\sigma > 1/2$, put for $u \in \mathrm{Lip}(I)$ and $x \in I$

$$\mathbf{H}_{s,i\tau}(u)(x) = \sum_{h \in \mathcal{H}} \mathrm{e}^{i\tau c(h)}|h'(x)|^s u(h(x)). \tag{13}$$

We have for the same $s$, $\tau$ and each $m \geq 1$, recalling the extension $c \colon \mathcal{H}^m \to \mathbb{R}$ from (9),

$$\mathbf{H}^m_{s,i\tau}(u)(x) = \sum_{h \in \mathcal{H}^m} \mathrm{e}^{i\tau c(h)} |h'(x)|^s u(h(x)).$$

Letting $f_1$ be the fixed point of $\mathbf{H}_{1,0}$ so that $\int_I f_1(x)\,\mathrm{d}x = 1$ (in fact, $f_1(x) = (\log 2)^{-1}(1+x)^{-1}$), we put for all $t \in \mathbb{R}$ and $\tau \in \mathbb{R}$

$$\widetilde{\mathbf{H}}_{1+it,i\tau}(u) = \frac{\mathbf{H}_{1+it,i\tau}(f_1 u)}{f_1}.$$

By definition, we have $\|\widetilde{\mathbf{H}}^m_{1+it,i\tau}(u)\|_{L^\infty} \leq \|u\|_{L^\infty}$, for all $m \in \mathbb{Z}^*_+$, and all real $t$ and $\tau$. It is not very difficult to show (see e.g. the proof of [1], Lemma 2) that for each $t \in \mathbb{R}$ there is $K(t) \geq 1$ so that

$$\mathrm{Lip}(\widetilde{\mathbf{H}}^m_{1+it,i\tau}(u)) \leq K(t)\|u\|_{L^\infty} + K\rho^m\,\mathrm{Lip}(u), \quad \forall m \in \mathbb{Z}^*_+, \forall \tau \in \mathbb{R}. \tag{14}$$

Inequalities such as the above are usually called Doeblin–Fortet (in probability) or Lasota–Yorke (in dynamics) inequalities. It will be convenient to use the following norm on $\mathrm{Lip}(I)$:

$$\|u\|_{\mathrm{Lip}} = \max\left\{\|u\|_{L^\infty}, \frac{1}{2\sup_{|t| \leq t_0} K(t)}\mathrm{Lip}(u)\right\}.$$

Indeed, recalling (12) and setting $n_0 = [\log K / \log(1/\rho)] + 1$, we have for all $\tau \in \mathbb{R}$ and all $t \in [-t_0, t_0]$

$$\begin{aligned}
\|\widetilde{\mathbf{H}}^m_{1+it,i\tau}\|_{\mathrm{Lip}} &\leq K\rho^m + 1, \quad \forall m \geq 1, \\
\|\widetilde{\mathbf{H}}^m_{1+it,i\tau}\|_{\mathrm{Lip}} &\leq 1, \qquad\qquad \forall m \geq n_0.
\end{aligned} \tag{15}$$

Finally, we give the definition of a diophantine cost:

**Definition.** *Let $\eta \geq 2$. The cost $c$ is diophantine of exponent $\eta$ for the finite subset $\mathcal{H}_0 \subset \mathcal{H}$ if there exists $\beta_0 \geq 1$ so that, for any sequences $\tau_k \in \mathbb{R}$, $t_k \in \mathbb{R}$, $\theta_k \in [0, 2\pi)$, with $\lim_{k\to\infty} |\tau_k| = \infty$ but $\sup_k |t_k| < \infty$, and for any $M \geq 1$ and $\beta \geq \beta_0$, there exist $k \geq 1$ and $x = h_x(x)$ for some $h_x \in \mathcal{H}^p_0$, with $p \geq 1$ minimal, so that*

$$\mathrm{dist}(\tau_k[\beta \log|\tau_k|]c(h_x) + t_k[\beta\log|\tau_k|]\log|h'_x| + p\theta_k, 2\pi\mathbb{Z}) \geq \frac{Mp}{|\tau_k|^\eta}.$$

It is easy to check that any diophantine cost is nonlattice.
The main result of this section is:

**Proposition 2.1.** *If the cost function $c$ is diophantine of exponent $\eta$ for $\mathcal{H}_0 \subset \mathcal{H}$, then, taking*

$$\alpha > \eta\left(2 + \frac{\log\sup_{\mathcal{H}_0}|h'|^{-1}}{\log(1/\rho)}\right)\left(1 + \frac{\log\sup_{\mathcal{H}_0}|h'|^{-1}}{\log(1/\rho)}\right),$$

*there exist $M_0 \geq 1$, and $\xi_0 \in (0, 1)$ so that for each $|\tau| \geq 2$, all $|t| \leq t_0$, and every $\sigma \geq 1 - \xi_0|\tau|^{-\alpha}$*

$$\|(\mathrm{Id} - \mathbf{H}_{\sigma+it,i\tau})^{-1}\|_{\mathrm{Lip}} \leq M_0|\tau|^\alpha.$$

(No growth assumption is required on $c$.)

Before we prove the proposition by a modification of the argument of Dolgopyat [7], as adapted by Melbourne [15], Section 3, to the case of infinitely many branches, we need further notation and a couple of preliminary lemmas. If $\mathcal{H}_0$ is a strict subset of $\mathcal{H}$, we let $I_0 = I_0(\mathcal{H}_0)$ be the invariant Cantor set for $T$ associated to $\mathcal{H}_0$, i.e., $I_0 = \{x \in I \mid T^m(x) \in \widetilde{I}_0, \forall m \in \mathbb{Z}_+\}$ for $\widetilde{I}_0 = \bigcup_{h \in \mathcal{H}_0} h([0,1])$. Then, we set $\mathrm{Lip}(I_0) =$ $\{u : I_0 \to \mathbb{C} \mid \mathrm{Lip}(u) < \infty\}$. Finally, for $\tau \in \mathbb{R}$ and $u \in L^\infty(I)$, we set for $x \in (0, 1]$, denoting by $h_x$ the element of $\mathcal{H}$ so that $x \in h_x([0, 1))$,

$$\mathcal{M}_{t,\tau} u(x) = |T'(x)|^{it} \mathrm{e}^{-i\tau c(h_x)} u(T(x)).$$

We may now state and prove the first lemma, which very roughly says that if the iterates of the transfer operator $\widetilde{\mathbf{H}}_{1+it,i\tau}$ decay too slowly, then the operator $\mathcal{M}_{t,\tau}$ has an approximate eigenfunction:

**Lemma 2.2.** *Let $I_0(\mathcal{H}_0)$ be associated to a finite set $\mathcal{H}_0 \subset \mathcal{H}$. Let $\eta_1 > 0$ and $\beta_0 \geq 1$. Then for each $\alpha_1 > \eta_1(2 + \frac{\log \sup_{\mathcal{H}_0} |h'|^{-1}}{\log(1/\rho)})$, there exist $\beta_1 > \beta_0$ and $K_0 \geq 1$ so that the following is true for each $|t| \leq t_0$ and every $|\tau| \geq 2$, setting $n(\tau) = [\beta_1 \log |\tau|]$:*

*Suppose that there exists $v_0 = v_{0,\tau,t} \in \mathrm{Lip}(I)$ with $\|v_0\|_{\mathrm{Lip}} \leq 1$ so that*

$$|\widetilde{\mathbf{H}}_{1+it,i\tau}^{jn(\tau)}(v_0)(x)| \geq 1 - \frac{1}{|\tau|^{\alpha_1}}, \quad \forall x \in I_0, j = 0, 1, 2. \tag{16}$$

*Then there exist $\theta_{\tau,t} \in [0, 2\pi)$ and $w_{\tau,t} : I_0 \to \mathbb{C}$, with $|w_{\tau,t}(x)| = 1$ for all $x$ and*

$$|\mathcal{M}_{t,\tau}^{n(\tau)} w_{\tau,t}(x) - \mathrm{e}^{i\theta_{\tau,t}} w_{\tau,t}(x)| \leq \frac{K_0}{|\tau|^{\eta_1}}, \quad \forall x \in I_0. \tag{17}$$

**Proof.** (See Lemma 3.12 and Section 3.3 in [15] adapted from [7], Section 8.) In this proof, we fix $t$ and $\tau$, and we write $n$ for $n(\tau)$. Letting $v_0$ be as in (16), put for $j = 0, 1, 2$

$$v_j = \widetilde{\mathbf{H}}_{1+it,i\tau}^{jn}(v_0), \qquad s_j = |v_j|.$$

Our assumption (16) implies that

$$1 - \frac{1}{|\tau|^{\alpha_1}} \leq s_j(x) \leq 1, \quad \forall x \in I_0, j = 0, 1, 2. \tag{18}$$

In particular, we may define $w_j(x) = \frac{v_j(x)}{s_j(x)}$ for $x \in I_0$ and $j = 0, 1, 2$ (and we have $|w_j| \equiv 1$ on $I_0$). Note for further use that (15) implies that there is a constant $K_1 \geq 1$ (which does not depend on $t$ or $\tau$) so that $\|w_j\|_{\mathrm{Lip}(I_0)} \leq K_1$ for $j = 0, 1, 2$ (first show that $\|v_j\|_{\mathrm{Lip}(I)}$ is uniformly bounded).

Since $s_1(x) = \frac{1}{w_1(x)} \widetilde{\mathbf{H}}_{1+it,i\tau}^n(s_0 w_0)(x)$, and $\widetilde{\mathbf{H}}_{1,0}(\mathbf{1}) = \mathbf{1}$ (with $\mathbf{1}$ the constant function $\equiv 1$), the bound (18) for $j = 1$ implies that for all $x \in I_0$

$$\sum_{h \in \mathcal{H}^n} \frac{f_1(h(x))}{f_1(x)} |h'(x)| \left(1 - \frac{|h'(x)|^{it}}{w_1(x)} \mathrm{e}^{i\tau c(h)} s_0(h(x)) w_0(h(x))\right) \leq \frac{1}{|\tau|^{\alpha_1}}.$$

The real part of each term in the above sum is nonnegative. Hence, using also that $f_1 \circ h/f_1$ is bounded from above and from below, uniformly in $h \in \mathcal{H}^n$, we can find a constant $K_2$ (which does not depend on $\tau$) so that for each $h \in \mathcal{H}^n$ and every $x \in I_0$

$$0 \leq 1 - s_0(h(x)) \Re\left(\frac{|h'(x)|^{it} \mathrm{e}^{i\tau c(h)} w_0(h(x))}{w_1(x)}\right) \leq \frac{K_2}{|h'(x)||\tau|^{\alpha_1}}.$$

Since $s_0(h(x)) \leq 1$, the above bound implies that for each $h \in \mathcal{H}^n$ and every $x \in I_0$

$$0 \leq 1 - \Re\left(\frac{|h'(x)|^{it} \mathrm{e}^{i\tau c(h)} w_0(h(x))}{w_1(x)}\right) \leq \frac{K_2}{|h'(x)||\tau|^{\alpha_1}}.$$

Using the fact that for any complex number $z$ of modulus 1 we have $|1 - z| = \sqrt{2}(1 - \Re z)^{1/2}$, we find a constant $K_3$, independent of $\tau$, so that for each $h \in \mathcal{H}^n$ and every $x \in I_0$

$$|w_1(x) - |h'(x)|^{it} e^{i\tau c(h)} w_0(h(x))| \leq \frac{K_3}{|h'(x)|^{1/2} |\tau|^{\alpha_1/2}}.$$

From now on, we restrict our attention to branches $h \in \mathcal{H}_0^n$. For such a branch, we have

$$|h'(x)| \geq K_4^{-n}, \tag{19}$$

where $K_4 = \sup_{h_0 \in \mathcal{H}_0} \sup |h_0'|^{-1} > 1$ depends only on $\mathcal{H}_0$. Therefore, recalling that $n = [\beta_1 \log |\tau|]$, if $\alpha_1$ and $\beta_1 > \beta_0$ satisfy

$$\alpha_1 - \beta_1 \log(K_4) > 2\eta_1, \tag{20}$$

then there is a constant $K_5$ (independent of $\tau$) so that for each $x \in I_0$, setting $h_x$ to be the element of $\mathcal{H}_0^n$ so that $x \in h_x([0,1))$,

$$|w_1(T^n(x)) - |(T^n)'(x)|^{-it} e^{i\tau c(h_x)} w_0(x)| \leq \frac{K_5}{|\tau|^{\eta_1}}. \tag{21}$$

A similar argument gives $K_6$, independent of $\tau$, so that if (20) holds and $x \in I_0$

$$|w_2(T^n(x)) - |(T^n)'(x)|^{-it} e^{i\tau c(h_x)} w_1(x)| \leq \frac{K_6}{|\tau|^{\eta_1}}. \tag{22}$$

Fix an arbitrary $x_0 \in I_0$ and define $\theta_0 = \theta_0(\tau, t)$ (recall that the $w_j$ depend on $\tau$ and $t$ and that $n = n(\tau)$) in $[0, 2\pi)$ by

$$e^{i\theta_0} = w_0(x_0)|(T^n)'(x_0)|^{-it} e^{i\tau c(h_{x_0})}.$$

Let $h_{x_0} \in \mathcal{H}_0^n$ be so that $x_0 \in h_{x_0}([0,1))$. Observe next that (21) and the fact that $T^n(x) = T^n(x_0)$ and $h_{x_0} = h_{x_{x_0}}$ for $x_{x_0} = h_{x_0}(T^n(x))$ imply that for all $x \in I_0$

$$|w_1(T^n(x)) - e^{i\theta_0}| \leq \frac{K_5}{|\tau|^{\eta_1}} + |w_0(x_{x_0}) - w_0(x_0)| + ||(T^n)'(x_{x_0})|^{-it} - |(T^n)'(x_0)|^{-it}|.$$

Now, on the one hand, since $|h''/h'| \leq \widehat{K}$, with $|h'| \leq K\rho^n$, and since $|t| \leq t_0$, and $|T^n(x) - T^n(x_0)| \leq 1$, there is a constant $K_7$, independent of $\tau$ and $h$, so that

$$||h_{x_0}'(T^n(x))|^{it} - |h_{x_0}'(T^n(x_0))|^{it}| \leq K_7 \cdot \rho^{n(\tau)},$$

and on the other hand, since $|x_{x_0} - x_0| \leq \rho^n$ and $x_{x_0} \in I_0$ if $x \in I_0$, we have

$$|w_0(x_{x_0}) - w_0(x_0)| \leq K_1 \cdot \rho^{n(\tau)}.$$

Therefore, if (20) holds, and in addition

$$\beta_1 > \frac{\eta_1}{\log(1/\rho)}, \tag{23}$$

then we have

$$|w_1(T^n(x)) - e^{i\theta_0}| \leq \frac{K_1 + K_5 + K_7}{|\tau|^{\eta_1}}, \quad \forall x \in I_0. \tag{24}$$

Next, defining $\theta_1 = \theta_1(t, \tau) \in [0, 2\pi)$ by

$$\mathrm{e}^{i\theta_1} = w_1(x_0)|(T^n)'(x_0)|^{-it}\mathrm{e}^{i\tau c(h_{x_0})},$$

the previous argument gives that if (20) and (23) hold then

$$|w_2(T^n(x)) - \mathrm{e}^{i\theta_1}| \le \frac{K_1 + K_6 + K_7}{|\tau|^{\eta_1}}, \quad \forall x \in I_0. \tag{25}$$

Putting together (24) and (25), we find for $\theta_{t,\tau} = \theta_0 - \theta_1$ and all $x \in I_0$

$$|\mathrm{e}^{-i\theta_{t,\tau}}w_1(T^n(x)) - w_2(T^n(x))| \le \frac{2K_1 + K_5 + K_6 + 2K_7}{|\tau|^{\eta_1}}. \tag{26}$$

Taking $\alpha_1$ and $\beta_1$ so that (20) and (23) hold, and substituting (26) into (22) we see that the function $w = w_1$ satisfies the conclusion of the lemma for $K_0 = 2K_1 + K_5 + 2K_6 + 2K_7$. $\qquad\square$

We need another lemma, which very roughly says that good decay for the iterates of $\widetilde{\mathbf{H}}_{1+it,i\tau}^{j_0}$ implies moderate $\tau$-growth for its resolvent:

**Lemma 2.3.** *Let $I_0(\mathcal{H}_0)$ be associated to a finite set $\mathcal{H}_0 \subset \mathcal{H}$. For each $\alpha_1 > 0$ every $\alpha > \alpha_1(1 + \frac{\log \sup_{\mathcal{H}_0} |h'|^{-1}}{\log(1/\rho)})$, and each $\beta_1 > 0$, there exists $M_0 \ge 1$ so that the following hold for each $|\tau| \ge 2$ and $|t| \le t_0$: Suppose that for each $v \in \mathrm{Lip}(I)$ with $\|v\|_{\mathrm{Lip}} \le 1$ there exists $x_0 \in I_0$ and $j_0 \le [3\beta_1 \log |\tau|]$ so that*

$$|\widetilde{\mathbf{H}}_{1+it,i\tau}^{j_0}(v)(x_0)| \le 1 - \frac{1}{|\tau|^{\alpha_1}}, \tag{27}$$

*then $\|(\mathrm{Id} - \mathbf{H}_{1+it,i\tau})^{-1}\|_{\mathrm{Lip}} \le M_0|\tau|^\alpha$.*

**Remark 2.4.** *The proof of Lemma 2.3 uses heavily the fact that $s = 1 + it$ and breaks down if $s = \sigma + it$ with $\sigma < 1$.*

**Proof of Lemma 2.3.** (See Lemma 3.13 in [15], adapted from [7], Section 7.) In this proof we fix $t$ and $\tau$ and write $n = n(\tau) = [3\beta_1 \log |\tau|]$. Let $v \in \mathrm{Lip}(I)$ be so that $\|v\|_{\mathrm{Lip}} \le 1$. It suffices to show that $(\mathrm{Id} - \widetilde{\mathbf{H}}_{1+it,i\tau})^{-1}(v)$ exists and $\|(\mathrm{Id} - \widetilde{\mathbf{H}}_{1+it,i\tau})^{-1}(v)\|_{\mathrm{Lip}} \le M_0|\tau|^\alpha$ for some $\alpha > 0$ and $M_0 \ge 1$ which do not depend on $\tau$ or $t$.

For $j_0 = j_0(v, \tau) \le n$ as in (27), we put

$$u_0 = \widetilde{\mathbf{H}}_{1+it,i\tau}^{j_0}(v) \quad \text{and} \quad u = \widetilde{\mathbf{H}}_{1+it,i\tau}^n(v).$$

We have $\|u_0\|_{L^\infty} \le 1$ and (recalling (15)) $\max\{\|u_0\|_{\mathrm{Lip}}, \|u\|_{\mathrm{Lip}}\} \le 1 + K\rho$. By (27), there is $x_0 \in I_0$ (depending on $v$) so that, putting $\bar{K} = 2(1 + K\rho)\sup_{|t| \le t_0} K(t)$

$$|u_0(x)| \le 1 - \frac{1}{2|\tau|^{\alpha_1}}, \quad \forall x \in I_{x_0} := \left\{ x \in I \ \Big| \ |x - x_0| \le \frac{1}{\bar{K}|\tau|^{\alpha_1}} \right\}. \tag{28}$$

Recall that $\mu_1$ is the absolutely continuous probability measure on $I$ with density $f_1$ (which is $T$-invariant). By definition, the dual of $\widetilde{\mathbf{H}}_{1,0}$ fixes $\mu_1$. Put $m_0 = m_0(\tau) = [\frac{\alpha_1 \log |\tau| + \log(K\bar{K})}{\log(1/\rho)}] + 1$, then the element $h_{x_0} \in \mathcal{H}_0^{m_0}$ so that $x_0 \in h_{x_0}([0, 1))$ is such that $h_{x_0}(I) \subset I_{x_0}$. Therefore, $\mu_1(I_{x_0}) \ge \mu_1(h_{x_0}(I))$. By definition of $\mu_1$, we have (recalling (12))

$$\mu_1(h_{x_0}(I)) \ge K_8^{-1}|h_{x_0}'(x_0)| \ge K_8^{-1}K_4^{-m_0(\tau)} \ge K_9^{-1}|\tau|^{-\alpha_1 \frac{\log(K_4)}{\log(1/\rho)}}, \tag{29}$$

with constants $K_8 \geq 1$, $K_9 \geq 1$, and $K_4 \geq 1$ (recall the proof of Lemma 2.2) independent of $\tau$ and $t$. Putting

$$\alpha_2 = \alpha_1 \left( \frac{\log(K_4)}{\log(1/\rho)} + 1 \right) > 0,$$

and decomposing $I = I_{x_0} \cup (I \setminus I_{x_0})$, we deduce from (28) and (29) that

$$\|u\|_{L^1(\mu_1)} \leq \|u_0\|_{L^1(\mu_1)} \leq \left( 1 - \frac{1}{2|\tau|^{\alpha_1}} \right) \mu_1(I_{x_0}) + 1 - \mu_1(I_{x_0}) = 1 - \frac{\mu_1(I_{x_0})}{2|\tau|^{\alpha_1}}$$

$$\leq 1 - K_9^{-1} |\tau|^{-\alpha_2}. \tag{30}$$

We next upgrade the $L^1$ estimate (30), first into an $L^\infty$ bound, and later into a Lipschitz estimate. For this, setting $n_1 = [\beta_2 \log |\tau|]$, for $\beta_2 > 1$ to be determined later, we get from the spectral decomposition (see e.g. [1] for references)

$$\widetilde{\mathbf{H}}_{1,0}^m(w) = \int w \, d\mu_1 + \mathbf{R}_{1,0}^m(w),$$

with $\|\mathbf{R}_{1,0}^m\|_{\mathrm{Lip}} \leq K_{10} \hat{\rho}^m$, for some $\hat{\rho} < 1$ and all $m \geq 1$, that

$$\|\widetilde{\mathbf{H}}_{1+it,i\tau}^{n_1(\tau)}(u)\|_{L^\infty} \leq \|\widetilde{\mathbf{H}}_{1,0}^{n_1(\tau)}|u|\|_{L^\infty} \leq \int |u| \, d\mu_1 + K_{10} \hat{\rho}^{n_1} \|u\|_{\mathrm{Lip}}$$

$$\leq 1 - K_9^{-1} |\tau|^{-\alpha_2} + (1 + K\rho) K_{10} \hat{\rho}^{n_1}.$$

Then, if $\beta_2 > \alpha_2$ is large enough (depending on $\rho$, $\hat{\rho}$, $K_9$ and $K_{10}$, but not on $\tau$) we have

$$\|\widetilde{\mathbf{H}}_{1+it,i\tau}^{n(\tau)+n_1(\tau)}(v)\|_{L^\infty} = \|\widetilde{\mathbf{H}}_{1+it,i\tau}^{n_1(\tau)}(u)\|_{L^\infty} \leq 1 - K_{11}^{-1} |\tau|^{-\alpha_2},$$

for $K_{11} \geq 1$ independent of $\tau$. Using (14), we get for $n_2(\tau) = [\beta_3 \log |\tau|]$ with large enough $\beta_3 > 3\beta_1 + \beta_2$ that

$$\|\widetilde{\mathbf{H}}_{1+it,i\tau}^{n_2(\tau)}(v)\|_{\mathrm{Lip}} \leq 1 - (2K_{11})^{-1} |\tau|^{-\alpha_2}.$$

Thus, since $v$ was arbitrary, $\|(\mathrm{Id} - \widetilde{\mathbf{H}}_{1+it,i\tau}^{n_2(\tau)})^{-1}\|_{\mathrm{Lip}} \leq 2K_{11} |\tau|^{\alpha_2}$. Finally, using

$$(\mathrm{Id} - \mathbf{A})^{-1} = (\mathrm{Id} + \mathbf{A} + \mathbf{A}^2 + \cdots + \mathbf{A}^{n_2 - 1})(\mathrm{Id} - \mathbf{A}^{n_2})^{-1},$$

and (15), we find for every $\alpha > \alpha_2$ a constant $M_0 \geq 1$, independent of $\tau$ and $t$, so that $\|(\mathrm{Id} - \widetilde{\mathbf{H}}_{1+it,i\tau})^{-1}\|_{\mathrm{Lip}} \leq M_0 |\tau|^\alpha$.    □

We may finally prove the proposition. The idea of the proof is that existence of approximate eigenfunctions contradicts the diophantine condition.

**Proof of Proposition 2.1.** The statement is trivial for $\sigma > 1$, because the spectral radius of $\mathbf{H}_{\sigma+it,i\tau}$ is then $< 1$ (see e.g. [1]). For $\sigma \leq 1$, we follow [15], Sections 3.2 and 3.3.

Let us first consider the case $\sigma = 1$, proceeding by contradiction. Let $I_0(\mathcal{H}_0)$ be associated to $\mathcal{H}_0 \subset \mathcal{H}$. Fix $\eta_1 > \eta$ and $\beta_0 \geq 1$. Then take $\alpha_1$, $\beta_1 > \beta_0$, and $K_0$ as in Lemma 2.2. Finally, take $\alpha$ and $M_0$ from Lemma 2.3. Assume that for each $M \geq M_0$ the bound $\|(\mathrm{Id} - \mathbf{H}_{1+it,i\tau})^{-1}\|_{\mathrm{Lip}} \leq M |\tau|^\alpha$ is violated for some $\tau = \tau(M)$ and $t = t(\tau)$ with $|\tau| > 2$ and $|t| \leq t_0$. By taking a sequence $M_k \to \infty$ we get sequences $t_k$ and $\tau_k$, with $|\tau_k|$ tending to infinity.

Then Lemma 2.3 implies that the hypothesis (16) of Lemma 2.2 is satisfied for each $(t, \tau) = (t_k, \tau_k)$ and for $\eta_1$. Therefore there are $\theta_k = \theta_{\tau_k, t_k} \in [0, 2\pi)$ and $w_k = w_{\tau_k, t_k} : I_0 \to \mathbb{C}$, with $|w_k(x)| = 1$ for all $x \in I_0$ and, setting $n_k = [\beta_1 \log |\tau_k|]$,

$$|\mathcal{M}_{t_k, \tau_k}^{n_k} w_k(x) - e^{i\theta_k} w_k(x)| \leq \frac{K_0}{|\tau_k|^{\eta_1}}, \quad \forall x \in I_0.$$

If $x = h_x(x) \in I_0$, for $h_x \in \mathcal{H}_0^p$, has minimal period $p \geq 1$, setting $c_x = c(h_x)$, we get

$$||h_x'|^{-it_k n_k} e^{-i\tau_k n_k c_x} w_k(x) - e^{ip\theta_k} w_k(x)| \leq \frac{pK_0}{|\tau_k|^{\eta_1}}.$$

Since $|w_k(x)| = 1$, we find integers $\ell_k(x) = \ell(t_k, \tau_k, x)$ with $|\ell_k| = \mathrm{O}(p|\tau_k| \log |\tau_k|)$ and a constant $D_0 \geq 1$, independent of $k$ and $x$, so that

$$|-t_k n_k \log |h_x'| - \tau_k n_k c_x - p\theta_k - 2\pi\ell_k(x)| \leq pD_0|\tau_k|^{-\eta_1}. \tag{31}$$

Since $\eta_1 > \eta$ and $\beta_0$ were arbitrary, this contradicts our diophantine assumption on $c$ when $k \to \infty$.

If $\sigma \in (1 - \xi_0|\tau|^{-\alpha_1}, 1)$, we put $s = \sigma + it$, and we write

$$(\mathrm{Id} - \widetilde{\mathbf{H}}_{s,i\tau})^{-1} = (\mathrm{Id} - \widetilde{\mathbf{H}}_{1+it,i\tau})^{-1}(\mathrm{Id} - \mathbf{A}_{s,\tau})^{-1}$$

with $\mathbf{A}_{s,\tau} = (\widetilde{\mathbf{H}}_{s,i\tau} - \widetilde{\mathbf{H}}_{1+it,i\tau})(\mathrm{Id} - \widetilde{\mathbf{H}}_{1+it,i\tau})^{-1}$. It is not very difficult to prove that for each $\sigma_0 > 1/2$ there is a constant $K_{12} \geq 1$ so that

$$\|\widetilde{\mathbf{H}}_{\sigma+it,i\tau} - \widetilde{\mathbf{H}}_{1+it,i\tau}\|_{\mathrm{Lip}} \leq K_{12}(1 - \sigma), \quad \forall \sigma \in (\sigma_0, 1], \forall |t| \leq t_0, \forall \tau \in \mathbb{R}.$$

(Use the bijection $\ell \mapsto h_\ell$ between $\ell \in \mathbb{Z}_+^*$ and $h_\ell \in \mathcal{H}$, from the introduction and the fact that $|h_\ell'| \leq \ell^{-2}$.) Thus, using the already treated case $\sigma = 1$, we get that $\|\mathbf{A}_{s,\tau}\|_{\mathrm{Lip}} \leq M_0 K_{12}(1 - \sigma)|\tau|^\alpha$ for all $\sigma > 1$ and all $|\tau| \geq 2$. This implies that there is $\xi_0 \in (0,1)$ so that $\|\mathbf{A}_{s,\tau}\|_{\mathrm{Lip}} \leq 1/2$ for $\sigma \geq 1 - \xi_0|\tau|^{-\alpha}$, and the result follows. $\qquad\square$

It remains to prove Lemma 1.5 stated in the Introduction:

**Proof of Lemma 1.5.** (See e.g. [7], Section 13 or [15], Corollary 2.4.) Consider the smallest $\mathcal{H}_0 \subset \mathcal{H}$ containing all points in the orbits of the four periodic points $x_j$, and let $\beta_0 \geq 1$ be fixed. Let $\eta > \eta_0$ and assume that $c$ is not diophantine of exponent $\eta$ for $\mathcal{H}_0$. It follows that there are $\beta \geq \beta_0$, $D \geq 1$, sequences $\tau_k$, $\theta_k$, and $t_k$, and integers $\ell_{k,j} = \ell_j(t_k, \tau_k)$ so that, putting $n_k = [\beta \log |\tau_k|]$,

$$|-t_k n_k \log |h_{i_j}'| - \tau_k n_k c(h_{i_j}) - p_j \theta_k - 2\pi\ell_{k,j}| < D|\tau_k|^{-\eta}, \quad j = 1, 2, 3, 4, \forall k. \tag{32}$$

Set $\tilde{\ell}_{k,j} = p_j \ell_{k,1} - \ell_{k,j} p_1 \in \mathbb{Z}$ for $j = 2, 3, 4$. We find by eliminating $\theta_k$ from (32) a constant $\widetilde{D} \geq 1$, independent of $k$, so that

$$|\tau_k n_k L_{1j} + t_k n_k \widehat{L}_{1j} + 2\pi\tilde{\ell}_{k,j}| < \widetilde{D}|\tau_k|^{-\eta}, \quad j = 2, 3, 4, \forall k. \tag{33}$$

If there is a sequence of $k \to \infty$ with $t_k = 0$, then we may assume that $L_{13} \neq 0$ and (up to taking a large enough $k$, i.e. a larger $\tau_k$) that $\tilde{\ell}_{k,3} \neq 0$. Then, since $\tilde{\ell}_{k,3} = \mathrm{O}(|\tau_k| \log |\tau_k|)$, eliminating $\tau_k n_k$ from (33) for $j = 2$ and $j = 3$ gives $\widehat{D}$ (independent of $k$) so that

$$\left|\tilde{\ell}_{k,3} \frac{L_{12}}{L_{13}} - \tilde{\ell}_{k,2}\right| < \widehat{D}|\tau_k|^{-\eta},$$

for all $k$ with $t_k = 0$, contradicting our diophantine assumption on $L_{12}/L_{13}$ when $k \to \infty$.

If this is not the case, we assume that $k$ is large enough so that all $t_k \neq 0$. Eliminating first $t_k n_k$ from (33), we get $\widehat{D}$ independent of $k$ such that

$$|\tau_k n_k(\widehat{L}_{1j}L_{12} - L_{1j}\widehat{L}_{12}) + 2\pi(\widehat{L}_{1j}\tilde{\ell}_{k,2} - \widehat{L}_{12}\tilde{\ell}_{k,j})| < \widehat{D}|\tau_k|^{-\eta}, \quad j = 3, 4, \forall k.$$

Eliminating $\tau_k n_k$, and dividing by $\widehat{L}_{12} \neq 0$ contradicts our strong diophantine assumption. (Using $\max(\tilde{\ell}_{k,2}, \tilde{\ell}_{k,3}) = \mathrm{O}(|\tau_k| \log |\tau_k|)$.) $\qquad\square$

## 3. Bounds on the moment generating function $\mathbb{E}_N(\mathrm{e}^{i\tau C})$ for smoothened costs

In this section, after recalling the methodology developed by Vallée in a series of papers (see e.g. [20]), as well as some useful lemmas from [1], we formulate in Corollary 3.5 a crucial consequence of Proposition 2.1.

The relevant sequence of Lévy moment generating functions is (see (46) and (47) in the next section) $\mathbb{E}_N(\mathrm{e}^{i\tau C})$ for $N \in \mathbb{Z}_+$, with $\tau \in \mathbb{R}$. To study this sequence of functions of $\tau$, we introduce a bivariate Dirichlet series:

$$S(s, i\tau) := \sum_{(p,q)\in\Omega} \frac{1}{q^s} \mathrm{e}^{i\tau C(p,q)}, \quad s = \sigma + it, \sigma > 2, t \in \mathbb{R}, \tau \in \mathbb{R}. \tag{34}$$

Now, on the one hand, denoting by $\mathbf{1}$ the constant function $\equiv 1$ on $I$, it is not very difficult to check (see e.g. [1]) that

$$S(2s, i\tau) = \mathbf{F}_{s,i\tau}(\mathrm{Id} - \mathbf{H}_{s,i\tau})^{-1}(\mathbf{1})(0), \tag{35}$$

with $\mathbf{F}_{s,i\tau}(u)(x) = \sum_{h \in \mathcal{H}'} \mathrm{e}^{i\tau c(h)} |h'(x)|^s u(h(x))$, where $\mathcal{H}' \subset \mathcal{H}$ contains all inverse branches of $T$ except $y \mapsto 1/(y+1)$.

On the other hand, the Perron formula of order two ([9], Theorem 2.7(b), see e.g. [1] for an application in the present context) gives that for each $D > 1$

$$\Psi_{i\tau}(N) = \frac{1}{2i\pi} \int_{D-i\infty}^{D+i\infty} S(2s, i\tau) \frac{N^{2s+1}}{s(2s+1)} \, \mathrm{d}s, \tag{36}$$

where $\Psi_{i\tau}(N)$ is the following Cesàro sum for the Dirichlet series:

$$\Psi_{i\tau}(N) = \sum_{M \leq N} \Phi_{i\tau}(M), \quad \text{with } \Phi_{i\tau}(M) = \sum_{(p,q)\in\Omega_M} \mathrm{e}^{i\tau C(p,q)}.$$

Our strategy will be to study $\Psi_{i\tau}(N)$, using (35) and (36) and spectral information on the transfer operators $\mathbf{H}_{s,i\tau}$. Clearly, $\mathbb{E}_N(\mathrm{e}^{i\tau C}) = \Phi_{i\tau}(N)/\Phi_0(N)$. However, exploiting directly estimates on $\Psi_{i\tau}(N)$ to get bounds on $\Phi_{i\tau}(N)$ seems difficult, and it is convenient to introduce instead auxiliary "smoothened" models as was done in [1], Section 4.2. The description given there was garbled – fortunately without consequences: All lemmas, propositions and theorems of [1] and [13] remain correct, except for [1], Lemma 14 and [13], Lemme 1(a), which should be replaced by Lemma 3.1. We give next the correct definition of the smoothening, as found by Eda Cesaratto (see [6]): for a function $\xi : \mathbb{Z}_+^* \to [0,1]$ and any integer $N \geq 1$, consider

$$\overline{\Omega_N}(\xi) = \bigcup_{N-[N\xi(N)] \leq Q \leq N} \Omega_Q \times \{Q\},$$

(noting that $\Omega_0 = \emptyset$), endowed with the uniform probability $\overline{\mathbb{P}}_N(\xi)$. Setting

$$\Pi(p, q, Q) = (p, q), \quad \text{for } (p, q, Q) \in \overline{\Omega_N}(\xi),$$

define

$$\overline{\Phi}_{i\tau}(\xi, N) = \overline{\Phi}_{i\tau}(N) := \sum_{(p,q,Q)\in\overline{\Omega}_N(\xi)} \mathrm{e}^{i\tau C(\Pi(p,q,Q))} = \sum_{Q=N-[N\xi(N)]}^{N} \sum_{q \leq Q} \sum_{(p,q)\in\Omega_q} \mathrm{e}^{i\tau C(p,q)}. \tag{37}$$

Then the moment generating function of the smoothened cost is just

$$\overline{\mathbb{E}}_N(\xi, \mathrm{e}^{i\tau C}) := \overline{\mathbb{E}}_N(\xi, \mathrm{e}^{i\tau C\circ\Pi}) = \frac{\overline{\Phi}_{i\tau}(N)}{\overline{\Phi}_0(N)}. \tag{38}$$

It is easy to check that

$$\overline{\Phi}_{i\tau}(N) = \Psi_{i\tau}(N) - \Psi_{i\tau}(N - [N\xi(N)]) - 1). \tag{39}$$

This implies that estimates on $\Psi_{i\tau}(N)$ give bounds on $\overline{\Phi}_{i\tau}(N)$. It remains to compare $\overline{\mathbb{P}}_N(\xi)$ with the primary object of interest, $\mathbb{P}_N$, and this is the purpose of Lemma 14 from [1], that we state here in its corrected form (see [6]):

**Lemma 3.1.** *There are $\widehat{M}_0 > 0$ and $\widehat{M} \geq 1$ so that for all $\xi$ with $(\xi(N))^{-1} \leq \widehat{M}_0 N / \log N$ we have*

$$|\overline{\mathbb{P}}_N(\xi)(\Pi^{-1}(E)) - \mathbb{P}_N(E)| \leq \widehat{M} \cdot \xi(N), \quad \forall E \subset \Omega_N, \forall N \in \mathbb{Z}_+^*.$$

**Sketch of proof.** First show that there is $\widetilde{M} \geq 1$ so that for all $N$ and all $(p, q) \in \Omega_{N-[N\xi(N)]}$

$$\left| \frac{\overline{\mathbb{P}}_N(\xi)(\Pi^{-1}((p,q)))}{\mathbb{P}_N((p,q))} - 1 \right| \leq \widetilde{M} \cdot \xi(N).$$

Then prove $\mathbb{P}_N(\Omega_N \setminus \Omega_{N-[N\xi(N)]}) = O(\xi(N))$ and $\overline{\mathbb{P}}_N(\xi)(\Pi^{-1}(\Omega_N \setminus \Omega_{N-[N\xi(N)]})) = O(\xi(N))$. $\square$

If $\xi$ satisfies the conditions of the above lemma, then it is easy to see that for any $F : \Omega_N \to \mathbb{C}$

$$|\overline{\mathbb{E}}_N(\xi, F \circ \Pi) - \mathbb{E}_N(F)| \leq \max |F| \cdot \widehat{M} \cdot \xi(N). \tag{40}$$

We shall work with two smoothenings. The first one, $\xi_1(N) = N^{-\gamma_0}$, for some small $\gamma_0 \in (0, 1)$ to be introduced in Lemmas 3.3 and 3.4, was used already in [1]. The second one appears only in the proof of Corollary 3.5 for $|\tau| \geq 2$.

In the proof of Theorem 1.9 in the next section, we will have to deal with $\overline{\mathbb{E}}_N(\xi_1, e^{i\tau C})$ for arbitrary $\tau \in \mathbb{R}$. The arguments for $\tau$ in a compact set (we shall use $|\tau| \leq 2$ to fix ideas) are the same as those in [1] for the case of lattice costs. Before stating the corresponding results from [1], let us mention an easy lemma which allows us to work with moment assumptions instead of moderate growth assumptions on the cost. If an operator has a simple eigenvalue $\lambda$ of modulus equal to its spectral radius, and if in addition the rest of the spectrum is contained in a disc of strictly smaller radius, we say that $\lambda$ is a *dominant eigenvalue*. For example, 1 is the dominant eigenvalue of $\mathbf{H}_{1,0}$ acting on Lip.

**Lemma 3.2.** *If $c$ has strong moments up to order $k$ for some $k \geq 1$, there exist $\nu_0 \in (0, 1)$ and $\nu_1 \in (0, 1/2)$ so that $(s, i\tau) \mapsto \mathbf{H}_{s,i\tau}$ is continuous (as an operator on* Lip*) on*

$$\mathcal{W} := \{(s, i\tau) \in \mathbb{C} \times i\mathbb{R} \mid |s - 1| < \nu_1, |\tau| < \nu_0\},$$

*and the dominant eigenvalue of $\mathbf{H}_{s,i\tau}$ acting on* Lip *is a continuous function $\lambda(s, i\tau)$ on $\mathcal{W}$, which is analytic in $s$ for each $\tau$. In addition, the (rank one) spectral projector $\mathbf{P}_{s,i\tau}$ for $\mathbf{H}_{s,i\tau}$ and $\lambda(s, i\tau)$ depends continuously on $(s, i\tau) \in \mathcal{W}$, and there is a uniquely defined continuous function $\sigma : (-\nu_0, \nu_0) \to \mathbb{C}$, with $\sigma(0) = 1$ and $\lambda(\sigma(i\tau), i\tau) \equiv 1$.*

*In fact, the functions $i\tau \mapsto \mathbf{H}_{s,i\tau}$, $i\tau \mapsto \mathbf{F}_{s,i\tau}$, $i\tau \mapsto \log \lambda(s, i\tau)$, $i\tau \mapsto \partial_s \lambda(s, i\tau)$, $i\tau \mapsto \mathbf{P}_{s,i\tau}$, and $i\tau \mapsto \sigma(i\tau)$ are $k$ times differentiable, uniformly in $s$ for $(s, i\tau) \in \mathcal{W}$. Their derivatives of order $0 \leq \ell \leq k$ are analytic functions of $s$, uniformly in each fixed $\tau$ for $(s, i\tau)$ in $\mathcal{W}$.*

**Proof.** See e.g. [1], Proposition 0, for the case of moderate growth, where all objects are analytic both in $s$ and $i\tau$. The continuous extension statements for $\lambda(s, i\tau)$ and $\mathbf{P}_{s,i\tau}$ follows from the (easily checked) fact that if we let $\mathbf{H}_{s,i\tau}$ act on the Banach space Lip$(I)$, then there is $\mathcal{W}$ so that $(s, i\tau) \mapsto \mathbf{H}_{s,i\tau}$ is continuous on $\mathcal{W}$ for the corresponding operator topology. Analyticity in $s$ is clear, and we have $\partial_s \lambda(1, 0) \neq 0$ as in [1], so that the implicit function theorem gives a continuous function $\sigma(i\tau)$ as claimed. Finally, up to taking smaller $\mathcal{W}$, the moment assumption on $c$ implies that, for each $|s - 1| < \nu_1$, and $\ell \leq k$, the $\ell$th derivative of $i\tau \mapsto \mathbf{H}_{s,i\tau}$, which is just $\sum_{h \in \mathcal{H}} (c(h))^{\ell} e^{i\tau c(h)} |h'(x)|^s u \circ h$, is a bounded operator on Lip, which is continuous in $(s, i\tau) \in \mathcal{W}$. $\qquad\square$

The following result is a small modification of Lemma 11 from [1]:

**Lemma 3.3.** *If the cost $c$ has strong moments up to order $k \geq 1$, letting $\nu_0 \in (0,1)$, $\mathbf{P}_{s,i\tau}$ and $\sigma : (-\nu_0, \nu_0) \to \mathbb{C}$ be as in Lemma 3.2, there exist $\tilde{\gamma}_0 \in (0,1)$ and $\gamma_1 \in (0,1/2)$ so that for each $\gamma_0 \in (0, \tilde{\gamma}_0)$, setting $\xi_1(N) = N^{-\gamma_0}$, we have for all $N \in \mathbb{Z}_+^*$*

$$\overline{\mathbb{E}}_N(\xi_1, e^{i\tau C}) = \frac{E(i\tau)}{E(0)\sigma(i\tau)} N^{2(\sigma(i\tau) - \sigma(0))}(1 + \mathrm{O}(N^{-\gamma_1})), \quad \forall |\tau| < \nu_0, \tag{41}$$

*where the $\mathrm{O}(N^{-\gamma_1})$ term is uniform in $\tau$, and $\tau \mapsto E(i\tau)$ is the $C^k$ function*

$$E(i\tau) = \frac{-1}{(\partial_s \lambda)(\sigma(i\tau), i\tau)} \mathbf{F}_{\sigma(i\tau), i\tau} \circ \mathbf{P}_{\sigma(i\tau), i\tau}(\mathbf{1})(0).$$

**Proof.** The correction in the definition of the smoothening corresponds to replacing the incorrect formula $[N\xi(N)]^{-1}(\Psi_{i\tau}(N) - \Psi_{i\tau}(N - [N\xi(N)]))$ for $\overline{\Phi}_{i\tau}(N)$ given in [1], (4.6), by (39). This is immaterial because the proof in [1] uses (38), so the factors $[N\xi(N)]$ cancel out, while the difference $\Psi_{i\tau}(N - [N\xi(N)]) - \Psi_{i\tau}(N - [N\xi(N)] - 1)$ is negligible. Since $i\tau$ is purely imaginary, we do not require the moderate growth assumption and we can apply Lemma 3.2 to adapt the proof of Lemma 11 in [1]. $\qquad\square$

The following claim is a small modification of Lemma 15 from [1], we shall apply it to $\nu_0$ from Lemma 3.3 and $L = 2$ in Section 4:

**Lemma 3.4.** *Let $c$ be nonlattice and let $\tilde{\gamma}_0 > 0$ be given by Lemma 3.3. For every $L > \nu_0 > 0$, there exist $\gamma_0 \in (0, \tilde{\gamma}_0)$, $\gamma_2 \in (0, 1/2)$ and $\widetilde{M} \geq 1$ so that, letting $\xi_1(N) = N^{-\gamma_0}$, we have for all $N \in \mathbb{Z}_+^*$*

$$|\overline{\mathbb{E}}_N(\xi_1, e^{i\tau C})| \leq \frac{\widetilde{M}}{N^{\gamma_2}}, \quad \forall |\tau| \in [\nu_0, L]. \tag{42}$$

**Proof.** By the last sentence of [1], Proposition 1, for a nonlattice cost, there exists $\sigma_1 < 1$ (depending on $\nu_0 > 0$ and $L$) so that $1 \notin \mathrm{sp}(\mathbf{H}_{\sigma+it,i\tau})$ (acting on Lip) for all $|\tau| \in (\nu_0, L)$, all $t \in \mathbb{R}$, and all $\sigma > \sigma_1$. Then, the proof of [1], Lemma 15, gives the claimed estimate. The moderate growth assumption is not used because $i\tau$ is purely imaginary. The correction in the definition of the smoothening is immaterial, as explained in Lemma 3.3. $\qquad\square$

The proofs of Lemmas 11 and 15 in [1] (implicit in Lemmas 3.3 and 3.4) use estimates on the growth of $(\mathrm{Id} - \mathbf{H}_{s,i\tau})^{-1}$, where $\tau$ is bounded, but where $s = \sigma + it$, with $|t|$ large, and $\sigma > \sigma_1$ with $\sigma_1 < 1$. These estimates are inspired from another important article of Dolgopyat [8], and use the fact that the Gauss map is "uniformly" away from a piecewise affine map (see the condition "UNI" in [1] for a precise formulation of this property and more details).

We now move on to "large" values of $\tau$. We shall prove that Proposition 2.1 implies the following estimate:

**Corollary 3.5.** *Let $c$ be diophantine, let $\alpha > 0$ be given by Proposition 2.1 and let $\xi_1(N) = N^{-\gamma_0}$ for $\gamma_0$ from Lemmas 3.3 and 3.4. For each $\alpha'' > \alpha' > \alpha$ there exists $K' \geq 1$ so that*

$$|\overline{\mathbb{E}}_N(\xi_1, e^{i\tau C})| \leq K' N^{-|\tau|^{-\alpha'}}, \quad \forall N \in \mathbb{Z}_+^*, \forall |\tau| \in [2, (\log N)^{1/\alpha''}]. \tag{43}$$

**Proof.** Fix $|\tau| \geq 2$ and $\alpha' > \alpha$, and introduce an auxiliary smoothening (only used in this proof)

$$\xi_2(N) = N^{-|\tau|^{-\alpha'}}. \tag{44}$$

By (40) and the triangle inequality we have for all $N \in \mathbb{Z}_+^*$

$$|\overline{\mathbb{E}}_N(\xi_1, e^{i\tau C}) - \overline{\mathbb{E}}_N(\xi_2, e^{i\tau C})| \leq \widehat{M}_1 N^{-\gamma_0} + \widehat{M}_2 N^{-|\tau|^{\alpha'}},$$

where $\widehat{M}_1$ and $\widehat{M}_2$ are uniform in $|\tau| \geq 2$. It thus suffices to prove the claimed estimate (43) for $\overline{\mathbb{E}}_N(\xi_2, e^{i\tau C})$.

As explained in the beginning of this section, our first goal is to obtain estimates on $\Psi_{i\tau}(N)$. Recall that we write $s = \sigma + it$ and let $\xi_0 \in (0, 1)$ be given by Proposition 2.1. We first claim that, up to taking a smaller $\xi_0$, for any $\mathcal{T} > 0$, and each $D > 1$ the function $s \mapsto S(2s, \tau)$ is holomorphic in

$$\mathcal{U}_{\mathcal{T}} := \{\sigma + it \mid |t| \leq \mathcal{T}, \sigma \in [1 - \xi_0|\tau|^{-\alpha}, D]\}.$$

Recalling (35), it suffices to study $(\mathrm{Id} - \mathbf{H}_{s,i\tau})^{-1}(\mathbf{1})$ for $(s, i\tau) \in \mathcal{U}_{\mathcal{T}}$. If $\xi_0$ is small enough, there is $t_0 > 0$, independent of $c$, so that if $|t| \geq t_0$, Theorem 2 in [1] gives $M \geq 1$ and $\bar{\alpha} \in (0, 1/5)$ (both independent of $t$, $\tau$, and $\sigma$) so that[7]

$$\sup|(\mathrm{Id} - \mathbf{H}_{s,i\tau})^{-1}(\mathbf{1})| \leq M|t|^{\bar{\alpha}}, \tag{45}$$

and in particular $s \mapsto \mathbf{F}_{s,i\tau}(\mathrm{Id} - \mathbf{H}_{s,i\tau})^{-1}(\mathbf{1})(0)$ is analytic in $\{s \in \mathcal{U}_{\mathcal{T}} \mid |t| \geq t_0\}$. If $|t| \leq t_0$, we may apply Proposition 2.1, using the diophantine condition, and we get that $s \mapsto \mathbf{F}_{s,i\tau}(\mathrm{Id} - \mathbf{H}_{s,i\tau})^{-1}(\mathbf{1})(0)$ is analytic in $\{s \in \mathcal{U}_{\mathcal{T}} \mid |t| \leq t_0\}$.

Next, by Cauchy's theorem

$$\int_{\partial \mathcal{U}_{\mathcal{T}}} S(2s, i\tau) \frac{N^{2s+1}}{s(2s+1)} \, \mathrm{d}s = 0, \quad \forall N \in \mathbb{Z}_+^*.$$

Clearly (45) implies that $\int_{\partial \mathcal{U}_{\mathcal{T}}, \Im s = \pm \mathcal{T}} S(2s, i\tau) \frac{N^{2s+1}}{s(2s+1)} \, \mathrm{d}s \to 0$ as $\mathcal{T} \to \infty$, uniformly in $N \in \mathbb{Z}_+^*$. By Perron's formula (36), the integral along the right-hand side border $\Re s = D$ of $\mathcal{U}_{\mathcal{T}}$ tends to $\Psi_{i\tau}(N)$. Finally, it is not very difficult to see that Proposition 2.1 implies that for each $\beta > \alpha$, there is $K_{13} \geq 1$ (depending on $\alpha$, $\beta$, and $\alpha''$) so that for all $N$ and all $|\tau| \in [2, (\log N)^{1/\alpha''}]$

$$\left| \int_{\partial \mathcal{U}_{\mathcal{T}}, \Re s = 1 - \xi_0|\tau|^{-\alpha}} S(2s, i\tau) \frac{N^{2s+1}}{s(2s+1)} \, \mathrm{d}s \right| \leq 4M_0|\tau|^\alpha N^{3-2|\tau|^{-\alpha}} \leq K_{13} N^{3-2|\tau|^{-\beta}}.$$

(First note that $|\tau|^{-\alpha} - |\tau|^{-\beta} \geq |\tau|^{-\alpha}(1 - 2^{\alpha-\beta})$, and then use that $\alpha'' > \alpha$ implies $\frac{\alpha}{\alpha''} \log(\log N) \leq 2(1 - 2^{\alpha-\beta}) \log N^{1-\alpha/\alpha''}$ for all large enough $N$, depending only on $\alpha < \beta < \alpha''$.) Combining the observations in this paragraph with (39) and the definition of $\xi_2(N)$ gives a constant $K_{14} \geq 1$, so that for all $N$ and all $|\tau| \in [2, (\log N)^{1/\alpha''}]$,

$$|\overline{\Phi}_{i\tau}(\xi_2, N)| \leq K_{14} N^{3-|\tau|^{-2\beta}}.$$

Now it is not very difficult to prove that there is $K_{15} > 0$ so that $\overline{\Phi}_0(\xi_2, N) \geq K_{15} N^{3-|\tau|^{-\alpha'}}$ for all $N$: we use the argument from [6]. First note that $\overline{\Phi}_0(\xi_2, N) = |\overline{\Omega}_N(\xi_2)|$. Since $|\Omega_Q| = 3Q^2(1 + \mathrm{O}(\log Q/Q))/\pi^2$ (see e.g. [1], Proof of Lemma 14), the definition of $\overline{\Omega}_N(\xi_2)$ together with the fact that $\xi_2(N)^{-1} = \mathrm{O}(N/\log N)$ give that $|\overline{\Omega}_N(\xi_2)| = (N^3 - (N - [N\xi_2(N)])]^3)(1 + \mathrm{O}(\xi_2(N)))/\pi^2$.

Taking $\alpha' > \beta > \alpha$, we end the proof by applying (38) for $\xi = \xi_2$. $\qquad\square$

---

[7] It is claimed in [1] that one may take $t_0 = \rho^{-2}$, but this is in fact not clear. Note in particular that the second and third inequalities in line 5 of [1], p. 362, are true only if $M_0$ is replaced by $M_0(M_1)^k$, so that $t_0$ must be taken large enough.

## 4. Proof of Theorem 1.9

By Lemma 3.1, it suffices to prove Theorem 1.9 for $\overline{\mathbb{P}}_N = \overline{\mathbb{P}}_N(\xi_1)$ and $\overline{\mathbb{E}}_N = \overline{\mathbb{E}}_N(\xi_1)$, where the smoothening is $\xi_1(N) = N^{-\gamma_0}$, with $\gamma_0 > 0$ as in Lemmas 3.3 and 3.4. Letting $\chi_J$ denote the characteristic function of the interval $J$, and recalling (6), our starting point to prove the first claim of the theorem is standard:

$$\overline{\mathbb{P}}_N((C(\cdot) - \mathbf{Q}(x, N)) \in J) = \sum_{(p,q,M) \in \overline{\Omega}_N} \chi_J(C(p, q) - \mathbf{Q}(x, N)) \cdot \overline{\mathbb{P}}_N(p, q, M)$$

$$= \sum_{(p,q,M) \in \Omega_N} \overline{\mathbb{P}}_N((p, q, M)) \cdot \frac{1}{2\pi} \int_{\mathbb{R}} e^{i\tau(C(\cdot) - \mathbf{Q}(x,N))} \hat{\chi}_J(\tau) \, d\tau$$

$$= \frac{1}{2\pi} \int_{\mathbb{R}} \hat{\chi}_J(\tau) e^{-i\tau \mathbf{Q}(x,N)} \overline{\mathbb{E}}_N(e^{i\tau C}) \, d\tau. \tag{46}$$

(We used $(2\pi)\psi(y) = \int \hat{\psi}(\tau) e^{i\tau y} \, d\tau$, where $\hat{\psi}(\tau) = \int e^{-i\tau x} \psi(x) \, dx$ is the Fourier transform of a locally supported bounded $\psi$.) Similarly, for the second claim, we shall use

$$\overline{\mathbb{E}}_N(\psi(C(\cdot) - \mathbf{Q}(x, N))) = \sum_{(p,q,Q) \in \overline{\Omega}_N} \psi(C(\cdot) - \mathbf{Q}(x, N)) \cdot \overline{\mathbb{P}}_N(p, q, Q)$$

$$= \frac{1}{2\pi} \int_{\mathbb{R}} \hat{\psi}(\tau) e^{-i\tau \mathbf{Q}(x,N)} \cdot \overline{\mathbb{E}}_N(e^{i\tau C}) \, d\tau. \tag{47}$$

Since $\hat{\chi}_J(\tau)$ does not decay fast enough when $\tau \to \infty$, it will be necessary to regularise $\chi_J$. For this reason, the proof of the second claim of Theorem 1.9 is easier, and we shall start with this.

### 4.1. The case of smooth $\psi$

By (47) it suffices to analyse

$$I(N) = \sqrt{\log N} \int_{\mathbb{R}} \hat{\psi}(\tau) e^{-i\tau \mathbf{Q}(x,N)} \cdot \overline{\mathbb{E}}_N(e^{i\tau C}) \, d\tau.$$

Recalling $\nu_0 \in (0, 1)$ from Lemma 3.3, let us decompose the real axis into

$$|\tau| < \nu_0, \qquad |\tau| \in [\nu_0, 2], \qquad |\tau| \in [2, L_N], \qquad |\tau| > L_N, \tag{48}$$

where $L_N > 2$ will be determined later. (We shall have $L_N \to \infty$ as $N \to \infty$.) The decomposition (48) gives rise to four integrals $I(N) = I_1(N) + I_2(N) + I_3(N) + I_4(N)$.

The integral $I_1(N)$ over $|\tau| < \nu_0$ is the dominant term, and can be handled by exploiting Lemma 3.3. More precisely, the saddle-point argument in [13], Section 3 (see [1], Section 5.1, for the lattice case), can be applied verbatim[8]: First note that the function $\mathcal{E}(i\tau) := E(i\tau)/(E(0)\sigma(i\tau))$ is $C^1$ (it is in fact $C^3$ by Lemma 3.2) on $(-\nu_0, \nu_0)$, so that $\mathcal{E}(i\tau) = 1 + \mathcal{O}(|\tau|)$. Note also that $\sigma(i\tau)$ is a $C^3$ function of $\tau$, by our assumption on the strong moments of $c$ and Lemma 3.2, with $\sigma(0) = 1$, $\sigma'(0) = \mu(c)/2$, and $\sigma''(0) = (\delta(c)^2)/2 \neq 0$ (see [1], Lemma 12, noting that we can replace $w$ by $i\tau$ in (4.13–4.14) there), and we can use a Taylor expansion of degree two with a remainder which is $\mathrm{O}(|\tau|^3)$. Then decompose $|\tau| < \nu_0$ into $|\tau| < \tau_N$ and $|\tau| \in [\tau_N, \nu_0)$ with

$$\tau_N = \left( \frac{\log \log N}{\delta_0 \log N} \right)^{1/2}, \tag{49}$$

---

[8]In particular, the moderate growth assumption on $c$ is not necessary in [13], Théorème 3. It is however necessary to assume that $c$ is nonlattice (this hypothesis is missing in [13]).

where $\delta_0 > 0$ depends on $\nu_0$, but not on $N$ or $\tau$ (see [1], (5.3)). (We may assume that $N$ is large enough to ensure $\tau_N < \nu_0$.) We find $M_{1,\psi} \geq 1$ (depending on $\sup_\tau |\hat{\psi}(\tau)| \leq \sup |\psi| \cdot |J|$, where $J = \text{supp}(\psi)$) so that for all $x \in \mathbb{R}$ and all $N \in \mathbb{Z}_+^*$

$$\left| \frac{I_1(N)}{2\pi} - \hat{\psi}(0) \frac{\mathrm{e}^{-x^2/2}}{\delta(c)\sqrt{2\pi}} \right| \leq \frac{M_{1,\psi}}{\sqrt{\log N}}.$$

(Of course, $\hat{\psi}(0) = \int \psi(y)\,\mathrm{d}y$.) Note that the error term $\mathrm{O}((\log N)^{-1/2})$ above cannot be improved. However, as observed in Remark 1.11, if $c$ enjoys strong moments of order 4 or higher, replacing $\mathbf{Q}(x, N)$ by a series with more terms is possible, and produces a smaller error term, giving an Edgeworth expansion.

The integral $I_2(N)$ over $|\tau| \in [\nu_0, 2]$ can be handled by exploiting Lemma 3.4, similarly to what was done previously in [13], Section 3 (see also [1], Section 5.2; for the case of lattice costs, and note that the moderate growth assumption on $c$ is not needed): It is easy to see that for each $\gamma_3 \in (0, \gamma_2)$, there exists $M_{2,\psi} \geq 1$ (depending on $\sup_\tau |\hat{\psi}(\tau)|$) so that for all $x$ and all $N$

$$|I_2(N)| = \sqrt{\log N} \left| \int_{|\tau| \in [\nu_0, 2]} \hat{\psi}(\tau) \mathrm{e}^{-i\tau \mathbf{Q}(x,N)} \overline{\mathbb{E}}_N(\mathrm{e}^{i\tau C})\,\mathrm{d}\tau \right| \leq M_{2,\psi} N^{-\gamma_3}.$$

Clearly, the error term above is $\mathrm{O}((\log N)^{-d})$ for arbitrarily large $d \geq 1/2$.

The last two integrals are more interesting. Let us assume that

$$r > \alpha + 1,$$

with $\alpha > 0$ from Proposition 2.1 (recall that $\alpha$ depends on $\eta$ from the diophantine condition). Letting $\alpha'' \in (\alpha, r-1)$, we put

$$L_N = (\log N)^{1/\alpha''}. \tag{50}$$

Then, Corollary 3.5 implies that for each $\alpha' \in (\alpha, \alpha'')$ there is $M_{3,\psi} \geq 1$, so that for all $x$ and all $N$

$$|I_3(N)| = \sqrt{\log N} \left| \int_{|\tau| \in [2, L_N]} \hat{\psi}(\tau) \mathrm{e}^{-i\tau \mathbf{Q}(x,N)} \overline{\mathbb{E}}_N(\mathrm{e}^{i\tau C})\,\mathrm{d}\tau \right|$$

$$\leq M_{2,\psi} \sqrt{\log N}\, L_N N^{-L_N^{-\alpha'}}$$

$$\leq M_{2,\psi} (\log N)^{1/2 + 1/\alpha''} \mathrm{e}^{-(\log N)^{1 - \alpha'/\alpha''}} \leq \frac{M_{3,\psi}}{(\log N)^{1/2}}. \tag{51}$$

The error term above is in fact $\mathrm{O}((\log N)^{-d})$ for arbitrarily large $d \geq 1/2$.

For each integer $m \leq r$ there is $M_\psi^{(m)}$ so that $|\hat{\psi}(\tau)| \leq M_\psi^{(m)} |\tau|^{-m}$ for all $|\tau| \geq 2$ (just use integration by parts). Finally since $r \geq \alpha' + 1$ by our assumption on $r$, and since $|\mathrm{e}^{-i\tau \mathbf{Q}(x,N)} \overline{\mathbb{E}}_N(\mathrm{e}^{i\tau C})| \leq 1$, we find, taking an integer $m \in (\alpha'' + 1, r]$, a constant $M_{4,\psi} \geq 1$ so that for all $N \in \mathbb{Z}_+^*$

$$|I_4(N)| = \sqrt{\log N} \left| \int_{|\tau| \geq L_N} \hat{\psi}(\tau) \mathrm{e}^{-i\tau \mathbf{Q}(x,N)} \overline{\mathbb{E}}_N(\mathrm{e}^{i\tau C})\,\mathrm{d}\tau \right|$$

$$\leq M_\psi^{(m)} \sqrt{\log N} (\log N)^{-(m-1)/\alpha''} \leq \frac{M_{4,\psi}}{(\log N)^{1/2}}. \tag{52}$$

Putting together the estimates on $I_1$, $I_2$, $I_3$ and $I_4$ ends the proof of the second claim of Theorem 1.9.

Up to taking a larger value of $m$ (if $r$ is large enough) the error term in (52) can be made $\mathrm{O}((\log N)^{-d})$ for arbitrarily large $d \geq 1/2$.

### 4.2. The case of $\chi_J$

We shall approximate the Dirac delta by using Gaussian distributions, writing, for small $\delta > 0$,

$$\Delta_\delta(x) = \delta^{-1} \Delta\left(\frac{x}{\delta}\right), \quad \text{with } \Delta(x) = \frac{1}{\sqrt{\pi}} e^{-x^2}.$$

For further use, note:

**Lemma 4.1.** *There is $D_0 \geq 1$ so that for every $\psi \in L^1(\mathbb{R})$ with $y\psi(y) \in L^1(\mathbb{R})$, setting $\psi_\delta = \psi * \Delta_\delta$, the Fourier transform of $\psi_\delta$ satisfies:*

$$|\hat{\psi}_\delta(\tau)| \leq D_0 \int |\psi(y)| \, dy \cdot e^{-\delta^2 \tau^2} \leq D_0 \int |\psi(y)| \, dy, \quad \forall \delta > 0, \tau \in \mathbb{R}, \tag{53}$$

$$\text{Lip}(\hat{\psi}_\delta) \leq D_0\left(\int |y\psi(y)| \, dy + \int |\psi(y)| \, dy\right), \quad \forall \delta > 0. \tag{54}$$

*In addition, if $\psi$ is Lipschitz, we have*

$$\sup_x |\psi(x) - \psi_\delta(x)| \leq D_0 \text{Lip}(\psi) \cdot \delta, \quad \forall \delta > 0. \tag{55}$$

**Proof.** To show (53) and (54), use $\hat{\psi}_\delta(\tau) = \hat{\psi}(\tau) \cdot \widehat{\Delta}_\delta(\tau)$, and recall that the Fourier transform $\widehat{\Delta}_\delta(\tau)$ of $\Delta_\delta$ is

$$\widehat{\Delta}_\delta(\tau) = e^{-\delta^2 \tau^2}.$$

To show (55), use $\int \Delta_\delta(y) \, dy = 1$ to write

$$|\psi(x) - \psi_\delta(x)| = \left|\int_{\mathbb{R}} \Delta_\delta(y)(\psi(x) - \psi(x - y)) \, dy\right| \leq \text{Lip}(\psi) \int_{\mathbb{R}} \Delta_\delta(y)|y| \, dy,$$

and note that $\int_{\mathbb{R}_+} \frac{y}{\sqrt{\pi}\delta} e^{-y^2/\delta^2} \, dy = O(\delta)$. (Just write $z = y/\delta$.) $\qquad \square$

Write $J = [a, b]$. For small $\delta \in (0, (b-a)^2/4)$, in view of applying the previous lemma, we first approximate $\chi_J$ by two compactly supported Lipschitz functions $\psi^+ = \psi_J^{+,\delta} : \mathbb{R} \to [0, 1]$ and $\psi^- = \psi_J^{-,\delta} : \mathbb{R} \to [0, 1]$, as follows. The function $\psi^+$ is $\equiv 1$ on the interval $[a - \sqrt{\delta}, b + \sqrt{\delta}]$, it is $\equiv 0$ outside of $[a - 2\sqrt{\delta}, b + 2\sqrt{\delta}]$, and it is affine with slope $\pm \delta^{-1/2}$ on the two remaining intervals. Similarly, $\psi^-$ is $\equiv 1$ on $[a + 2\sqrt{\delta}, b - 2\sqrt{\delta}]$ it is $\equiv 0$ outside of $[a + \sqrt{\delta}, b - \sqrt{\delta}]$, and it is affine with slope $\pm \delta^{-1/2}$ on the two remaining intervals. We have that $\int \psi^\pm(y) \, dy = |J| + O(\sqrt{\delta})$ and $\int |y|\psi^\pm(y) \, dy \leq 4|J|^2$. In addition

$$\overline{\mathbb{E}}_N(\psi^-(C(\cdot) - \mathbf{Q}(x, N))) \leq \overline{\mathbb{P}}_N((C(\cdot) - \mathbf{Q}(x, N)) \in J) \leq \overline{\mathbb{E}}_N(\psi^+(C(\cdot) - \mathbf{Q}(x, N))).$$

Next, we consider the regularisation by convolution

$$\psi_\delta^\pm = \psi_J^{\pm,\delta} * \Delta_\delta, \quad \text{with } \delta = \delta_N = (\log N)^{-2\varepsilon},$$

with $\varepsilon > 0$ to be determined later. Since $\text{Lip}(\psi^\pm) = \delta^{-1/2}$, the bound (55) from Lemma 4.1 gives

$$|\overline{\mathbb{E}}_N(\psi^\pm(C(\cdot) - \mathbf{Q}(x, N))) - \overline{\mathbb{E}}_N(\psi_\delta^\pm(C(\cdot) - \mathbf{Q}(x, N)))| \leq D_0\sqrt{\delta} \leq D_0(\log N)^{-\varepsilon}.$$

Therefore, by (46), it suffices to analyse

$$I_\delta^\pm(N) = \sqrt{\log N} \int_{\mathbb{R}} \hat{\psi}_\delta^\pm(\tau) e^{-i\tau \mathbf{Q}(x, N)} \overline{\mathbb{E}}_N(e^{i\tau C}) \, d\tau.$$

We consider the decomposition (48) of $\mathbb{R}$, for $L_N$ to be determined later, and the four associated integrals $I_\delta^\pm(N) = I_{\delta,1}^\pm(N) + I_{\delta,2}^\pm(N) + I_{\delta,3}^\pm(N) + I_{\delta,4}^\pm(N)$.

Introducing $\tau_N$ like in (49), and using in addition the uniform bound (54) on the Lipschitz constant of $\hat{\psi}_\delta^\pm(\tau)$ in order to see that

$$\sup_{|\tau| \le \tau_N} |\hat{\psi}_\delta^\pm(0) - \hat{\psi}_\delta^\pm(\tau)| \le M_J \cdot \tau_N,$$

we find $M_{1,J} \ge 1$ (using Lemma 3.2 as in Section 4.1 and the weak claim in (53) to bound $\sup |\hat{\psi}^\pm|$) so that for all $x \in \mathbb{R}$ and all $N \in \mathbb{Z}_+^*$

$$\left| \frac{I_{\delta,1}^\pm(N)}{2\pi} - \hat{\psi}_\delta^\pm(0) \frac{e^{-x^2/2}}{\delta(c)\sqrt{2\pi}} \right| \le \frac{M_{1,J}}{(\log(N))^{1/2}}. \tag{56}$$

(The term $\mathrm{O}((\log N)^{-1/2})$ in the above expression cannot be improved.) Note that

$$\hat{\psi}_\delta^\pm(0) = \int \psi_\delta^\pm(y)\,\mathrm{d}y = \int \psi_\delta^\pm(y)\,\mathrm{d}y = |J| + \mathrm{O}(\sqrt{\delta_N}) = |J| + \mathrm{O}((\log N)^{-\varepsilon}). \tag{57}$$

Just like in Section 4.1, since $c$ is nonlattice, for each $\gamma_3 \in (0, \gamma_2)$, there exists $M_{2,J} \ge 1$ (using the weak bound from (53)) so that for all $x \in \mathbb{R}$ and all $N \in \mathbb{Z}_+^*$

$$|I_{\delta,2}^\pm(N)| \le M_{2,J} N^{-\gamma_3}. \tag{58}$$

Take $\alpha'' > \alpha$, with $\alpha$ from Corollary 3.5, and put $L_N = (\log N)^{1/\alpha''}$. (We may assume $\alpha'' > 2$.) Then, Corollary 3.5 (with the weak bound from (53)) implies that (see (51)) for each $\alpha' \in (\alpha, \alpha'')$ there is $M_{3,J} \ge 1$ so that for all $x \in \mathbb{R}$ and all $N \in \mathbb{Z}_+^*$

$$|I_{\delta,3}^\pm(N)| = \sqrt{\log N} \left| \int_{|\tau| \in [2, L_N]} \hat{\psi}_\delta^\pm(\tau) e^{-i\tau \mathbf{Q}(x,N)} \overline{\mathbb{E}}_N(e^{i\tau C})\,\mathrm{d}\tau \right| \le \frac{M_{3,J}}{(\log N)^{1/2}}. \tag{59}$$

Finally if

$$2\varepsilon < (\alpha'')^{-1}$$

then[9] $\delta_N L_N > \log \log N$ and, since $|e^{-i\tau \mathbf{Q}(x,N)} \overline{\mathbb{E}}_N(e^{i\tau C})| \le 1$ we find, using the strong bound from (53), that there are constants $\widetilde{M_{4,J}} \ge 1$ and $M_{4,J} \ge 1$ so that for all $x \in \mathbb{R}$ and all $N \in \mathbb{Z}_+^*$

$$|I_{\delta,4}^\pm(N)| = \sqrt{\log N} \left| \int_{|\tau| \ge L_N} \hat{\psi}_\delta^\pm(\tau) e^{-i\tau \mathbf{Q}(x,N)} \overline{\mathbb{E}}_N(e^{i\tau C})\,\mathrm{d}\tau \right|$$

$$\le \widetilde{M_{4,J}} \sqrt{\log N} L_N e^{-(\delta L_N)^2} \le \widetilde{M_{4,J}} \sqrt{\log N} L_N e^{-(\log \log N)^2} \le \frac{M_{4,J}}{(\log N)^{1/2}}. \tag{60}$$

Putting together (56), (58)–(60) we have proved the first claim of Theorem 1.9, for $\varepsilon \in (0, 1/(2\alpha''))$.

## Appendix. The centered and odd Euclidean algorithms

Let us describe two variants $\mathcal{K}$ and $\mathcal{O}$ of the classical continued fraction (1), for rational $x \in (0, 1]$. Both of them are of the form

$$x = \frac{1}{m_1 + \varepsilon_1/(m_2 + \varepsilon_2/(\cdots + \varepsilon_{P-1}/m_P))}, \quad m_j \in \mathbb{Z}_+^*, P \in \mathbb{Z}_+^*, \varepsilon_J \in \{-1, 1\}. \tag{61}$$

---

[9]We may take a smaller value of $\delta_N$, but our argument requires $\inf_N(\delta_N L_N) > 0$.

The first one is the *centered* algorithm $\mathcal{K}$, for which all digits $m_j$ are $\geq 2$. It is described by the following centered division algorithm, for integers $p$ and $q$ with $1 \leq p \leq q/2$: write $q = mp + \varepsilon r$, with $m \in \mathbb{Z}_+^*$, $\varepsilon \in \{-1, 1\}$, and integer $r$ with $\varepsilon r \in [-p/2, p/2]$. If $r = 0$ we take $m_{1,\mathcal{K}}(p/q) = m$, $P_{\mathcal{K}}(p/q) = 1$, and we are done. Otherwise, we put $m_1 = m$, $\varepsilon_1 = \varepsilon$, $p_2 = r_1 = r$, and $q_2 = p_1$, and we iterate until $r_{P_{\mathcal{K}}} = 0$, constructing the $m_{j,\mathcal{K}}(p/q)$ and $\varepsilon_{j,\mathcal{K}}(p/q)$ along the way. The associated dynamical system on $(0, 1/2]$ satisfies $T_{\mathcal{K}}(p/q) = r/p$ and is just

$$T_{\mathcal{K}}(x) = \left| \frac{1}{x} - A_{\mathcal{K}}\left(\frac{1}{x}\right) \right|,$$

where $A_{\mathcal{K}}(y)$ is the nearest integer to $y$, i.e., the unique integer $m$ so that $y - m \in [-1/2, 1/2)$.

The second one is the *odd* algorithm $\mathcal{O}$, for which all digits $m_j$ are odd. It is described by the following odd division algorithm, for integers $1 \leq p \leq q$: write $q = mp + \varepsilon r$, with odd $m \in \mathbb{Z}_+^*$, $\varepsilon \in \{-1, 1\}$, and an integer $r$ with $\varepsilon r \in [-p, p]$. If $r = 0$ we take $m_{1,\mathcal{O}}(p/q) = m$, $P_{\mathcal{O}}(p/q) = 1$, and we are done. Otherwise, we put $m_1 = m$, $\varepsilon_1 = \varepsilon$, $p_2 = r_1 = r$, and $q_2 = p_1$, and we iterate until $r_{P_{\mathcal{O}}} = 0$, constructing the $m_{j,\mathcal{O}}(p/q)$ and $\varepsilon_{j,\mathcal{O}}(p/q)$ along the way. The associated dynamical system on $(0, 1]$ satisfies $T_{\mathcal{O}}(p/q) = r/p$ and is just

$$T_{\mathcal{O}}(x) = \left| \frac{1}{x} - A_{\mathcal{O}}\left(\frac{1}{x}\right) \right|,$$

where $A_{\mathcal{O}}(y)$ is the nearest odd integer to $y$, i.e., the unique odd integer $m$ so that $y - m \in [-1, 1)$.

We refer e.g. to [1], Section 2, for more information on these two algorithms and their associated interval maps $T_{\mathcal{K}}$ and $T_{\mathcal{O}}$. Letting $\mathcal{H}_{\mathcal{K}}$ and $\mathcal{H}_{\mathcal{O}}$ denote the set of inverse branches of $T_{\mathcal{K}}$ and $T_{\mathcal{O}}$, respectively, it turns out that the corresponding transfer operators $\mathbf{H}_{s,i\tau,\mathcal{K}}$ and $\mathbf{H}_{s,i\tau,\mathcal{O}}$ enjoy the same properties as those of the operator $\mathbf{H}_{s,i\tau}$ associated to the ordinary Euclidean division and (1). (See [1], Section 2, for the invariant densities $f_{1,\mathcal{O}}$ and $f_{1,\mathcal{K}}$, the constants $\rho_{\mathcal{O}}$ and $\rho_{\mathcal{K}}$, etc., and also [1], Proposition 0, as well as [1], Section 3.4, for the condition "UNI.") In particular, all statements in Sections 2 and 3 of the present paper hold, replacing the Gauss map $T$ by $T_{\mathcal{O}}$ or $T_{\mathcal{K}}$, up to changing the constants. The proof in Section 4 can be followed for both algorithms, and finally we see that Theorem 1.9 is true also for the centered and odd algorithms, up to replacing $\mu(c)$ and $\delta(c)$ by appropriate constants $\mu_{\mathcal{K}}(c)$, $\delta_{\mathcal{K}}(c) > 0$ or $\mu_{\mathcal{O}}(c)$, $\delta_{\mathcal{O}}(c) > 0$.

## Acknowledgments

Part of this work was done during the trimester "Time at Work" in Institut Henri Poincaré, Paris, 2005. Warm thanks to F. Naud, who pointed out the paper of I. Melbourne, to E. Cesaratto who found a mistake in [1] and explained how to fix it, to E. Breuillard, S. Gouëzel, F. Ledrappier, and B. Vallée for useful comments. A. H. benefitted from a Gauthier–Villars Fellowship, via AFFDU. Both authors are partially supported by the ACI-DynamicAL, CNRS.